\documentclass{article}

\usepackage{latexsym}
\usepackage{amssymb}
\usepackage{amsthm}
\usepackage{amssymb,amsmath}
\usepackage{MnSymbol}
\usepackage{mathtools}
\usepackage{stmaryrd}
\usepackage{skak}
\usepackage{pigpen}

\newcommand{\bpr}{\noindent{\em Proof\/. }}
\newcommand{\epr}{\hspace*{\fill}$\meddiamond$\medskip}

\newcommand{\bprthone}{\noindent{\em Proof of Theorem~\ref{mainone}\/. }}

\newcommand{\bprthtwo}{\noindent{\em Proof of Theorem~\ref{maintwo}\/. }}

\numberwithin{equation}{section}

\newtheorem{thm}{Theorem}
\newtheorem{lem}{Lemma}
\newtheorem{cor}{Corollary}
\newtheorem{prop}{Proposition}

\newtheorem{conj}{Conjecture}

\theoremstyle{definition}
\newtheorem{rem}{Remark}

\title{Prime bound of a graph}

\author{Abderrahim Boussa\"{\i}ri\thanks{Facult\'e des Sciences A\"{\i}n Chock, 
D\'epartement de Math\'ematiques et Informatique, Km 8 route d'El Jadida, 
BP 5366 Maarif, Casablanca, Maroc; {\tt aboussairi@hotmail.com}.}\and 
Pierre Ille\thanks{Institut de Math\'ematiques de Luminy, 
CNRS -- UMR 6206, 163 avenue de Luminy, Case~907, 13288 Marseille Cedex 09, France; {\tt ille@iml.univ-mrs.fr}.}
\thanks{Centre de recherches math\'ematiques, Universit\'e de Montr\'eal, Case postale 6128, 
Succursale Centre-ville, Montr\'eal, Qu\'ebec, Canada H3C 3J7.}}

\begin{document}

\maketitle

\begin{abstract}
Given a graph $G$, a subset $M$ of $V(G)$ is a module of $G$ if for each 
$v\in V(G)\setminus M$, $v$ is adjacent to all the elements of $M$ or to none of them. 
For instance, $V(G)$, $\emptyset$ and $\{v\}$ ($v\in V(G)$) are modules of $G$ called trivial. 
Given a graph $G$, $m(G)$ denotes the largest integer 
$m$ such that there is a module $M$ of $G$ which is a clique or a stable set in $G$ with  
$|M|=m$. 
A graph $G$ is prime if $|V(G)|\geq 4$ and if all its modules are trivial. 
The prime bound of $G$ is the smallest integer $p(G)$ such that there is a prime graph $H$ with $V(H)\supseteq V(G)$, 
$H[V(G)]=G$ and $|V(H)\setminus V(G)|=p(G)$. 
We establish the following. 
For every graph $G$ such that $m(G)\geq 2$ and $\log_2(m(G))$ is not an integer, 
$p(G)=\lceil\log_2(m(G))\rceil$. 
Then, we prove that for every graph $G$ such that $m(G)=2^k$ where $k\geq 1$, 
$p(G)=k$ or $k+1$. 
Moreover $p(G)=k+1$ if and only if $G$ or its complement admits $2^k$ isolated vertices. 
Lastly, we show that $p(G)=1$ for every non-prime graph $G$ such that $|V(G)|\geq 4$ and $m(G)=1$. 
\end{abstract}

\medskip


\medskip


\section{Introduction}

A {\em graph} $G=(V(G),E(G))$ is constituted by a {\em vertex set} $V(G)$ and an {\em edge set} $E(G)\subseteq\binom{V(G)}{2}$. 
Given a set $S$, $K_S=(S,\binom{S}{2})$ is the {\em complete} graph on $S$ whereas $(S,\emptyset)$ is the {\em empty} graph. 
Let $G$ be a graph. 
With each $W\subseteq V(G)$ associate the {\em subgraph} $G[W]=(W,\binom{W}{2}\cap E(G))$ of $G$ induced by $W$. 
A graph $H$ is an {\em extension} of $G$ if $V(H)\supseteq V(G)$ and $H[V(G)]=G$. 
Given $p\geq 0$, a $p$-extension of $G$ is an extension $H$ of $G$ such that $|V(H)\setminus V(G)|=p$. 
The {\em complement} of $G$ is the graph $\overline{G}=(V(G),\binom{V(G)}{2}\setminus E(G))$. 
A subset $W$ of $V(G)$ is a {\em clique} (respectively a {\em stable set}) in $G$ if $G[W]$ is complete (respectively empty). 
The largest cardinality of a clique (respectively a stable set) in $G$ is the {\em clique number} (respectively the {\em stability number}) of $G$, denoted by 
$\omega(G)$ (respectively $\alpha(G)$). 
Given $v\in V(G)$, the {\em neighbourhood} $N_G(v)$ of $v$ in $G$ is the family $\{w\in V(G):\{v,w\}\in E(G)\}$.
Its {\em degree} is $d_G(v)=\left|N_G(v)\right|$. 
We will consider $N_G$ as a function from $V(G)$ in $2^{V(G)}$. 
A vertex $v$ of $G$ is {\em isolated} if $N_G(v)=\emptyset$. 
The family of isolated vertices of $G$ is denoted by ${\rm Iso}(G)$. 

We use the following notation. 
Let $G$ be a graph. 
For $v\neq w\in V(G)$, 
\begin{equation*}
(v,w)_G=
\begin{cases}
0&\text{if $\{v,w\}\not\in E(G)$},\\
1&\text{if $\{v,w\}\in E(G)$}. 
\end{cases}
\end{equation*}
Given $W\subsetneq V(G)$, $v\in V(G)\setminus W$ and $i\in\{0,1\}$, $(v,W)_G=i$ means $(v,w)_G=i$ for every $w\in W$. 
Given $W,W'\subsetneq V(G)$, with $W\cap W'=\emptyset$, and $i\in\{0,1\}$, $(W,W')_G=i$ means $(w,W')_G=i$ for every $w\in W$. 
Given $W\subsetneq V(G)$ and $v\in V(G)\setminus W$, $v\sim_GW$ means that there is $i\in\{0,1\}$ such that $(v,W)_G=i$. 
The negation is denoted by $v\not\sim_GW$.

Given a graph $G$, a subset $M$ of $V(G)$ is a {\em module} of $G$ if for each 
$v\in V(G)\setminus M$, we have $v\sim_GM$. 
For instance, $V(G)$, $\emptyset$ and $\{v\}$ ($v\in V(G)$) are modules of $G$ called {\em trivial}. 
Clearly, if $|V(G)|\leq 2$, then all the modules of $G$ are trivial. 
On the other hand, if $|V(G)|=3$, then $G$ admits a non-trivial module. 
A graph $G$ is then said to be {\em prime} if $|V(G)|\geq 4$ and if all its modules are trivial. 
For instance, given $n\geq 4$,  the {\em path} $(\{1,\ldots,n\},\{\{p,q\}:|p-q|=1\})$ is prime. 
Given a graph $G$, $G$ and $\overline{G}$ share the same modules. 
Thus $G$ is prime if and only if $\overline{G}$ is. 

Let $S$ be a set with $|S|\geq 2$. 
Given $p\geq 1$, consider a $p$-extension $G$ of $\overline{K_S}$. 
If $|S|\geq 2^{p}$, then $G$ is not prime. 
Indeed, for each $s\in S$, we have $N_G(s)\subseteq V(G)\setminus S$ because $S$ is a stable set in $G$. 
So if $|S|> 2^{p}$, then $(N_G)_{\restriction S}:S\longrightarrow 2^{V(G)\setminus S}$ is not injective. 
Thus $\{s,t\}$ is a non-trivial module of $G$ for $s\neq t\in S$ such that $N_G(s)=N_G(t)$. 
Furthermore, if $(N_G)_{\restriction S}:S\longrightarrow 2^{V(G)\setminus S}$ is injective and if $|S|=2^{p}$, then there is 
$s\in S$ such that $N_G(s)=\emptyset$. 
Therefore $s\in{\rm Iso}(G)$ and $V(G)\setminus\{s\}$ is a non-trivial module of $G$. 
On the other hand, the following is well known and is easily verified 
(see Sumner~\cite[Theorem~2.45]{S71} and also Corollary~\ref{c2Xstable} below). 
Given a set $S$ with $|S|\geq 2$, 
$\overline{K_S}$ admits a prime $\lceil\log_2(|S|+1)\rceil$-extension. 
This is extended to any graph by Brignall \cite[Theorem 3.7]{B07} as follows. 

\begin{thm}\label{tBrignall}
A graph $G$, with $|V(G)|\geq 2$, admits a prime extension $H$ 
$$\text{such that}\ |V(H)\setminus V(G)|\leq\lceil\log_2(|V(G)|+1)\rceil.$$ 
\end{thm}

Following Theorem~\ref{tBrignall}, we introduce the notion of prime bound. 
Let $G$ be a graph. 
The {\em prime bound} of $G$ is the smallest integer $p(G)$ such that $G$ admits a prime 
$p(G)$-extension. 
Obviously $p(G)=0$ when $G$ is prime. 

A prime extension $H$ of $G$ is {\em minimal} 
\cite{O90,G97,Z03,BHV04,GO07} if for every $W\subsetneq V(H)$ such that $H[W]$ is prime, $H[W]$ does not admit an induced subgraph isomorphic to $G$. 
Given a graph $G$, a prime $p(G)$-extension of $G$ is clearly minimal. 
By Theorem~\ref{tBrignall}, $p(G)\leq\lceil\log_2(|V(G)|+1)\rceil$. 
Observe also that $p(G)=p(\overline{G})$ for every graph $G$. 
By considering the clique number and the stability number, Brignall \cite[Conjecture 3.8]{B07} conjectured the following. 

\begin{conj}\label{Brignall}
For each graph $G$ with $|V(G)|\geq 2$, 
$$p(G)\leq\lceil\log_2(\max(\omega(G),\alpha(G))+1)\rceil.$$
\end{conj}

We answer the conjecture positively by refining the notions of clique number and of stability number as follows. 
Given a graph $G$, the {\em modular clique number} of $G$ is the largest integer 
$\omega_M(G)$ such that there is a module $M$ of $G$ which is a clique in $G$ with  
$|M|=\omega_M(G)$. 
The {\em modular stability number} of $G$ is 
$\alpha_M(G)=\omega_M(\overline{G})$. 
Obviously $\omega_M(G)\leq \omega(G)$ and $\alpha_M(G)\leq\alpha(G)$. 
For convenience, set 
$$m(G)=\max(\alpha_M(G),\omega_M(G)).$$ 
We establish 

\begin{thm}\label{mainone}
For every  graph $G$ such that $m(G)\geq 2$, 
$$\lceil\log_2(m(G))\rceil\leq p(G)\leq\lceil\log_2(m(G)+1)\rceil.$$
\end{thm}

On the one hand, it follows that $p(G)=\lceil\log_2(m(G))\rceil$ for a graph $G$ such that $m(G)\geq 2$ and 
$\log_2(m(G))$ is not an integer. 
On the other, if $\log_2(m(G))$ is a positive integer, that is, $m(G)=2^k$ where $k\geq 1$, then 
$p(G)=k$ or $k+1$. 
We prove 

\begin{thm}\label{maintwo}
For every graph $G$ such that $m(G)=2^k$ where $k\geq 1$, 
\begin{equation*}
\text{$p(G)=k+1$ if and only if $|{\rm Iso}(G)|=2^k$ or $|{\rm Iso}(\overline{G})|=2^k$.} 
\end{equation*}
\end{thm}

Lastly, we show that $p(G)=1$ for every non-prime graph $G$ such that $|V(G)|\geq 4$ and $m(G)=1$ (see Proposition~\ref{m(G)=1}). 

The case of directed graphs is quite different. 
Recall that a {\em tournament} $T=(V(T),A(T))$ is a directed graph such that 
$|\{(v,w),(w,v)\}\cap A(T)|=1$ for any $v\neq w\in V(T)$. 
For instance, the {\em 3-cycle} $C_3=(\{1,2,3\},\{(1,2),(2,3),(3,1)\})$ is a tournament. 
A tournament $T$ is {\em transitive} provided that for any $u,v,w\in V(T)$, if $(u,v),(v,w)\in A(T)$, then $(u,w)\in A(T)$. 
Given a tournament $T$, a subset $I$ of $V(T)$ is an {\em interval} \cite{ST93,I97} (or a {\em clan} \cite{EHR99}) of $T$ if for any $x,y\in I$ and $v\in V(T)\setminus I$, we have: $(x,v)\in A(T)$ if and only if $(y,v)\in A(T)$. 
Once again, $V(T)$, $\emptyset$ and $\{v\}$ ($v\in V(T)$) are intervals of $T$ called trivial. 
For instance, all the intervals of $C_3$ are trivial. 
On the other hand, a transitive tournament with at least 3 vertices admits a non-trivial interval. 
A tournament with at least 3 vertices is {\em indecomposable} \cite{ST93,I97} ( or {\em simple} 
\cite{EFHM72,EHM72,M72}) if all its intervals are trivial. 
The {\em indecomposable bound} of a tournament $T$ is defined as the prime bound of a graph. It is still denoted by $p(T)$. 

\begin{thm}[\cite{EFHM72,EHM72,M72}]
For a tournament $T$ with $|V(T)|\geq 3$, $$p(T)\leq 2.$$ 
Moreover
$p(T)=2$ if and only if 
$T$ is transitive and $|V(T)|$ is odd.
\end{thm}

\section{Preliminaries}

We begin with the well known properties of the modules of a graph (for example, see \cite[Lemma~3.4, Theorem~3.2, Lemma~3.9]{EHR99}). 

\begin{prop}\label{properties}
Let $G$ be a graph. 
\begin{enumerate}
\item Given $W\subseteq V(G)$, if $M$ is a module of $G$, then $M\cap W$ is a module of $G[W]$. 
\item Given a module $M$ of $G$, if $N$ is a module of $G[M]$, then $N$ is a module of $G$. 
\item If $M$ and $N$ are modules of $G$, then $M\cap N$ is a module of $G$.
\item If $M$ and $N$ are modules of $G$ such that $M\cap N\neq\emptyset$, then $M\cup N$ is a module of $G$.
\item If $M$ and $N$ are modules of $G$ such that $M\setminus N\neq\emptyset$, then $N\setminus M$ is a module of $G$.
\item If $M$ and $N$ are modules of $G$ such that $M\cap N=\emptyset$, then there is $i\in\{0,1\}$ such that $(M,N)_G=i$.
\end{enumerate}
\end{prop}

Given a graph $G$, a partition $P$ of $V(G)$ is a {\em modular partition} of $P$ if each element of $P$ is a module of $G$. 
Let $P$ be such a partition. 
Given $M\neq N\in P$, there is $i\in\{0,1\}$ such that $(M,N)_G=i$ by the last assertion of Proposition~\ref{properties}. 
This justifies the following definition. 
The {\em quotient} of $G$ by $P$ is the graph $G/P$ defined on $V(G/P)=P$ by 
$(M,N)_{G/P}=(M,N)_G$ for $M\neq N\in P$. 
We use the following properties of the quotient (for example, see 
\cite[Theorems~4.1--4.3, Lemma~4.1]{EHR99}). 

\begin{prop}\label{pquotient}
Given a graph $G$, consider a modular partition $P$ of $G$.
\begin{enumerate}
\item Given $W\subseteq V(G)$, if $|W\cap X|=1$ for each $X\in P$, then $G[W]$ and 
$G/P$ are isomorphic. 
\item For any module $M$ of $G$, $\{X\in P:M\cap X\neq\emptyset\}$ is a module of $G/P$. 
\item For any module $Q$ of $G/P$, the union $\cup Q$ of the elements of $Q$ is a module of $G$. 
\end{enumerate}
\end{prop}

Given a graph $G$, with each non-empty module $M$ of $G$, associate the modular partition 
$P_M=\{M\}\cup\{\{v\}:v\in V(G)\setminus M\}$. 
Given $m\in M$, the corresponding quotient $G/P_M$ is isomorphic to 
$G[(V(G)\setminus M)\cup\{m\}]$ by the first assertion of 
Proposition~\ref{pquotient}. 
To associate a unique quotient with any graph and to characterize the corresponding quotient, the following strengthening of the notion of module is introduced. 
Given a graph $G$, a module $M$ of $G$ is said to be {\em strong} provided that for every module $N$ of $G$, we have: if $M\cap N\neq\emptyset$, then 
$M\subseteq N$ or $N\subseteq M$. 
We recall the following well known properties of the strong modules of a 
graph (for example, see 
\cite[Theorem~3.3]{EHR99}).

\begin{prop}\label{strong}
Let $G$ be a graph.
\begin{enumerate}
\item $V(G)$, $\emptyset$ and $\{v\}$, $v\in V$, are strong modules of $G$.
\item For a strong module $M$ of $G$ and for $N\subseteq M$, $N$ is a strong module of $G$ if and only if $N$ is a strong module of $G[M]$. 
\end{enumerate}
\end{prop}

With each graph $G$, we associate the family $\Pi(G)$ of the maximal strong modules of $G$ under inclusion among the proper and non-empty strong modules of $G$. 
The modular decomposition theorem is stated as follows. 

\begin{thm}[Gallai \cite{G67,MP01}]\label{Gallai}
For a graph $G$ with $|V(G)|\geq 2$, the family $\Pi(G)$ realizes a modular partition of $G$. 
Moreover, the corresponding quotient $G/\Pi(G)$ is complete, empty or prime. 
\end{thm}

Given a graph $G$ with $|V(G)|\geq 2$, we denote by $\mathbb{S}(G)$ the family of the non-empty strong modules of $G$. 
As a direct consequence of the definition of a strong module, we obtain that 
the family $\mathbb{S}(G)$ endowed with inclusion, denoted by 
$(\mathbb{S}(G),\subseteq)$, is a tree called the {\em modular decomposition tree}  \cite{CM05} of $G$. 
Given $M\in\mathbb{S}(G)$ with $|M|\geq 2$, it follows from 
Proposition~\ref{strong} that $\Pi(G[M])\subseteq\mathbb{S}(G)$. 
Furthermore, given $W\subseteq V(G)$ 
(respectively $W\subsetneq V(G)$), 
$(\{M\in\mathbb{S}(G):M\supseteq W\},\subseteq)$ (respectively 
$(\{M\in\mathbb{S}(G):M\supsetneq W\},\subseteq)$) is a total order. 
Its smallest element is denoted by $W\!\!\uparrow$ (respectively 
$W\!\!\twoheaduparrow$). 
By Proposition~\ref{strong}, if $M\in\mathbb{S}(G)\setminus\{V(G)\}$, then 
$M\in\Pi(G[M\!\!\twoheaduparrow])$. 

Let $G$ be a graph with $|V(G)|\geq 2$. 
Using Theorem~\ref{Gallai}, we label $\mathbb{S}(G)\setminus\{\{v\}:v\in V(G)\}$ by the function $\lambda_G$ defined as follows. 
For each $M\in\mathbb{S}(G)$ with $|M|\geq 2$, 
\begin{equation*}
\lambda_G(M)=
\begin{cases}
\filledmedsquare\ \ \text{if $G[M]/\Pi(G[M])$ is complete,}\\
\medsquare\ \ \text{if $G[M]/\Pi(G[M])$ is empty,}\\
\sqcup\ \ \text{if $G[M]/\Pi(G[M])$ is prime.}
\end{cases}
\end{equation*}

\begin{figure}[!h]
\begin{center}
\setlength{\unitlength}{0.7cm}
\begin{picture}(14,6)

\put(0,3){$\bullet$}\put(-0.5,2.8){$s_1$}
\put(2,3){$\bullet$}\put(1.5,2.8){$s_2$}
\put(4,3){$\bullet$}\put(3.5,2.8){$s_3$}
\put(6,3){$\bullet$}\put(6.3,2.8){$s_4$}
\put(8,3){$\bullet$}\put(8.3,2.8){$s_5$}

\put(4.1,3){\oval(9.6,2)}
\put(0.1,2.2){$S$}

\put(4.1,0.1){\line(-4,3){4}}
\put(4.1,0.1){\line(-2,3){2}}
\put(4.1,0.1){\line(0,1){3}}
\put(4.1,0.1){\line(2,3){2}}
\put(4.1,0.1){\line(4,3){4}}

\put(12,3){$\bullet$}\put(11.6,2.8){$c_1$}
\put(14,3){$\bullet$}\put(14.3,2.8){$c_2$}
\put(13,1.5){$\bullet$}\put(12.87,2.1){$c_3$}

\put(13.1,2.5){\oval(4,2.5)}
\put(14.3,1.6){$C$}

\put(13.1,0.1){\line(0,1){1.5}}
\put(13.1,0.1){\line(-1,3){1}}
\put(13.1,0.1){\line(1,3){1}}

\put(13.1,1.6){\line(-2,3){1}}
\put(13.1,1.6){\line(2,3){1}}
\put(12.1,3.1){\line(1,0){2}}

\put(4,0){$\bullet$}\put(3.6,-0.2){$a$}
\put(13,0){$\bullet$}\put(13.4,-0.2){$b$}
\put(4.1,0.1){\line(1,0){9}}

\end{picture}
\caption{
$\mathbb{S}(G)\setminus\{\{v\}:v\in V(G)\}=\{C,S,V(G)\}$, 
$\lambda_G(C)=\filledmedsquare$, $\lambda_G(S)=\medsquare$ and 
$\lambda_G(V(G))=\sqcup$.
\label{fig1}
}
\end{center}
\end{figure}
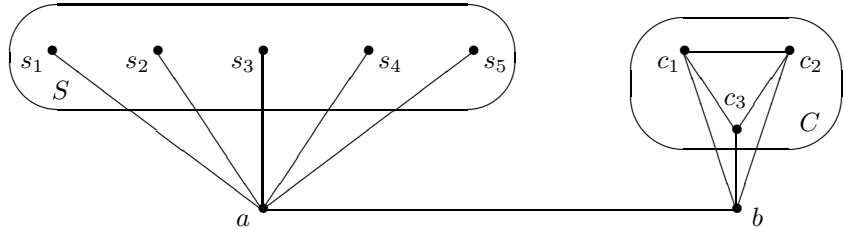

In Figure~\ref{fig1}, we depict a graph $G$ defined on 
$V(G)=\{a,b\}\cup\{c_1,c_2,c_3\}\cup\{s_1,\ldots,s_5\}$ by
\begin{equation*}
\begin{cases}
\text{$C=\{c_1,c_2,c_3\}$ is a clique in $G$,}\\
\text{$S=\{s_1,\ldots,s_5\}$ is a stable set in $G$,}\\
(a,S)_G=(a,b)_G=(b,C)_G=1,\\
(a,C)_G=(C,S)_G=(b,S)_G=0.
\end{cases}
\end{equation*} 
For $i\in\{1,2,3\}$ and $j\in\{1,\ldots,5\}$, $G[\{a,b,c_i,s_j\}]$ is a path and hence is prime. 
Thus the non-trivial modules of $G$ are the subsets $M$ of $C$ or of $S$ with $|M|\geq 2$. 
It follows that 
$$\mathbb{S}(G)=\{\{v\}:v\in V(G)\}\cup\{C,S,V(G)\}$$ 
\begin{equation*}
\text{and}\ 
\begin{cases}
\lambda_G(C)=\filledmedsquare,\\
\lambda_G(S)=\medsquare,\\
\lambda_G(V(G))=\sqcup.
\end{cases}
\end{equation*}

\section{Modular clique number and modular stability number}

Given a non-prime graph $G$, we are looking for a prime extension $H$ of $G$. 
We have to break the non-trivial modules of $G$ as modules of $H$. 
Precisely, given a non-trivial module $M$ of $G$, there must exist $v\in V(H)\setminus V(G)$ such that $v\not\sim_HM$. 
Moreover, we have only to consider the minimal non-trivial modules of $G$ under inclusion. 

Given a graph $G$, denote by $\mathcal{M}(G)$ the family of modules $M$ of $G$ such that 
$\left|M\right|\geq 2$, and denote by $\mathcal{M}_{{\rm min}}(G)$ the family of the minimal elements of $\mathcal{M}(G)$ under inclusion. 

\begin{lem}\label{min}
Let $G$ be a graph. 
For every $M\in\mathcal{M}(G)$, $M\in\mathcal{M}_{{\rm min}}(G)$ if and only if either $\left|M\right|=2$ or $G[M]$ is prime. 
\end{lem}

\bpr
To begin, consider $M\in\mathcal{M}_{{\rm min}}(G)$. 
It follows from the second assertion of Proposition~\ref{properties} that all the modules of $G[M]$ are trivial. 
Thus $G[M]$ is prime if $\left|M\right|\geq 4$. 
Assume that $\left|M\right|\leq 3$. 
Since a graph on 3 vertices admits a non-trivial module, we obtain $\left|M\right|=2$. 

Conversely, consider $M\in\mathcal{M}(G)$ such that $\left|M\right|=2$ or $G[M]$ is prime. 
Clearly $M\in\mathcal{M}_{{\rm min}}(G)$ when $\left|M\right|=2$. 
So assume that $G[M]$ is prime and consider $N\in\mathcal{M}(G)$ such that $N\subseteq M$. 
By the first assertion of Proposition~\ref{properties}, $N$ is a module of $G[M]$. 
As $G[M]$ is prime, $N=M$ and hence $M\in\mathcal{M}_{{\rm min}}(G)$. 
\epr

Let $G$ be a graph. Following Lemma~\ref{min}, we denote by $\mathcal{P}(G)$ the family of modules $M$ of $G$ such that $G[M]$ is prime. 
In Figure~\ref{fig1}, $\mathcal{P}(G)=\emptyset$. 

\begin{lem}\label{prime}
Let $G$ be a graph. 
For any $P\in\mathcal{P}(G)$ and $M\in\mathcal{M}(G)$, either $P\cap M=\emptyset$ or 
$P\subseteq M$. 
\end{lem}

\bpr
Assume that $P\cap M\neq\emptyset$. 
By the first assertion of Proposition~\ref{properties}, $P\cap M$ is a module of $G[P]$. 
Since $G[P]$ is prime, either $\left|P\cap M\right|=1$ or $P\cap M=P$. 
Suppose for a contradiction that $\left|P\cap M\right|=1$. 
As $\left|M\right|\geq 2$, $M\setminus P\neq\emptyset$. 
By the second to last assertion of Proposition~\ref{properties}, $P\setminus M$ is a module of $G$. 
By the first assertion, $P\setminus M$ would be a non-trivial module of $G[P]$. 
It follows that $P\cap M=P$. 
\epr 

Given a graph $G$, consider $u\neq v\in V(G)$ and $v\neq w\in V(G))$ such that there are 
$M_{\{u,v\}},M_{\{v,w\}}\in\mathcal{M}_{{\rm min}}(G)$ with $u,v\in M_{\{u,v\}}$ and 
$v,w\in M_{\{v,w\}}$. 
First, assume that $M_{\{u,v\}}\in\mathcal{P}(G)$. 
By Lemma~\ref{prime}, $M_{\{u,v\}}\subseteq M_{\{v,w\}}$ and hence $M_{\{u,v\}}=M_{\{v,w\}}$ 
because $M_{\{v,w\}}\in\mathcal{M}_{{\rm min}}(G)$. 
Similarly $M_{\{u,v\}}=M_{\{v,w\}}$ when $M_{\{v,w\}}\in\mathcal{P}(G)$. 
Second, assume that $\left|M_{\{u,v\}}\right|=\left|M_{\{v,w\}}\right|=2$, that is, 
$M_{\{u,v\}}=\{u,v\}$ and $M_{\{v,w\}}=\{v,w\}$. 
By interchanging $G$ and $\overline{G}$, assume that $\{u,v\}$ is a stable set in $G$. 
We obtain $(u,\{v,w\})_G=0$ and hence $(w,\{u,v\})_G=0$. 
Thus $\{u,v,w\}$ is a stable set in $G$.  
Furthermore $\{u,v,w\}$ is a module of $G$ by the fourth assertion of Proposition~\ref{properties}. 
Since $\{u,v,w\}$ is a stable set in $G$, $\{u,w\}$ is a module of $G[\{u,v,w\}]$. 
By the second assertion of Proposition~\ref{properties}, 
$\{u,w\}$ is a module of $G$.  
By Lemma~\ref{min}, $\{u,w\}\in\mathcal{M}_{{\rm min}}(G)$. 
In both cases, there exists $M_{\{u,w\}}\in\mathcal{M}_{{\rm min}}(G)$ such that 
$u,w\in M_{\{u,w\}}$. 
Consequently, the binary relation $\approx_G$ defined on $V(G)$ by 
\begin{equation*}
u\approx_Gv \text{\ if\ }
\begin{cases}
u=v\\
\text{or}\\
\text{$u\neq v$ and there is $M_{\{u,v\}}\in\mathcal{M}_{{\rm min}}(G)$ with 
$u,v\in M_{\{u,v\}}$},
\end{cases}
\end{equation*} 
for any $u,v\in V(G)$, is an equivalence relation. 
The family of the equivalence classes of $\approx_G$ is denoted by 
$\mathfrak{M}(G)$. 
In Figure~\ref{fig1}, $\mathfrak{M}(G)=\{C,S,\{a\},\{b\}\}$. 

\begin{lem}\label{partition}
Let $G$ be a graph. 
For every $C\in\mathfrak{M}(G)$ such that $\left|C\right|\geq 2$, one and only one of the following holds
\begin{itemize}
	\item $C\in\mathcal{P}(G)$,
	\item $C$ is a maximal module of $G$ under inclusion among the modules of $G$ which are cliques in $G$,
	\item $C$ is a maximal module of $G$ under inclusion among the modules of $G$ which are stable sets in $G$.
\end{itemize}
\end{lem}

\bpr
Consider $C\in\mathfrak{M}(G)$ such that $\left|C\right|\geq 2$. 
Clearly $C$ satisfies at most one of the three assertions above. 
For any $u\neq v\in C$, there is $M_{\{u,v\}}\in\mathcal{M}_{{\rm min}}(G)$ such that 
$u,v\in M_{\{u,v\}}$. 
By Lemma~\ref{min}, either $\left|M_{\{u,v\}}\right|=2$ or $M_{\{u,v\}}\in\mathcal{P}(G)$. 

First, assume that there are $u\neq v\in C$ such that $M_{\{u,v\}}\in\mathcal{P}(G)$. 
Let $w\in C\setminus\{u,v\}$. 
By Lemma~\ref{prime}, $M_{\{u,v\}}\subseteq M_{\{u,w\}}$ and hence 
$M_{\{u,v\}}=M_{\{u,w\}}$ 
because $M_{\{u,w\}}\in\mathcal{M}_{{\rm min}}(G)$. 
Thus $C=M_{\{u,v\}}\in\mathcal{P}(G)$. 

Second, assume that $M_{\{u,v\}}=\{u,v\}$ for any $u\neq v\in C$. 
Given $u\neq v\in C$, assume that $\{u,v\}$ is a stable set in $G$ by interchanging $G$ and $\overline{G}$. 
As observed above, $\{u,v,w\}$ is a module of $G$ and a stable set in $G$ for every 
$w\in C\setminus\{u,v\}$. 
Therefore $C$ is a stable set in $G$ and it follows from the fourth assertion of Proposition~\ref{properties} that $C$ is a module of $G$. 
Furthermore consider $u\in C$ and $x\in V(G)\setminus C$. 
Since $u\not\approx_Gx$, $\{u,x\}$ is not a module of $G$ and hence there is $y\in V(G)\setminus\{u,x\}$ such that $y\not\sim_G\{u,x\}$. 
If $y\not\in C$, then $y\not\sim_GC\cup\{x\}$ and $C\cup\{x\}$ is no longer a module of $G$. 
If $y\in C$, then $(y,x)_G=1$ because $(y,u)_G=0$. 
Thus $C\cup\{x\}$ is no longer a stable set in $G$. 
Consequently $C$ is a maximal module of $G$ among the modules of $G$ which are stable sets in $G$.
\epr

Let $G$ be a graph. 
Following Lemma~\ref{partition}, denote by $\mathcal{C}(G)$ the family of the maximal elements of 
$\mathcal{M}(G)$ under inclusion among the elements of 
$\mathcal{M}(G)$ which are cliques in $G$, and denote by $\mathcal{S}(G)$ the family of the maximal elements of 
$\mathcal{M}(G)$ under inclusion among the elements of 
$\mathcal{M}(G)$ which are stable sets in $G$. 
In Figure~\ref{fig1}, $\mathcal{C}(G)=\{C\}$ and $\mathcal{S}(G)=\{S\}$. 
The next is a simple consequence of Lemma~\ref{partition}. 

\begin{cor}\label{cpartition}
For a graph $G$, $\{C\in\mathfrak{M}(G):|C|\geq 2\}=\mathcal{C}(G)\cup\mathcal{S}(G)\cup\mathcal{P}(G)$.
\end{cor}

\bpr
By Lemma~\ref{partition}, 
$$\{C\in\mathfrak{M}(G):|C|\geq 2\}\subseteq\mathcal{C}(G)\cup\mathcal{S}(G)\cup\mathcal{P}(G).$$
For the opposite inclusion, consider $C\in\mathcal{C}(G)\cup\mathcal{S}(G)\cup\mathcal{P}(G)$. 
First, assume that $C\in\mathcal{C}(G)\cup\mathcal{S}(G)$. 
By interchanging $G$ and $\overline{G}$, assume that $G$ is a clique of $G$. 
For any $c\neq d\in C$, $\{c,d\}$ is a module of $G[C]$. 
By the second assertion of Proposition~\ref{properties}, $\{c,d\}$ is a module of $G$. 
Thus $\{c,d\}\in\mathcal{M}_{{\rm min}}(G)$ and hence $c\approx_G d$. 
Therefore there is $D\in\mathfrak{M}(G)$ such that $D\supseteq C$. 
By Lemma~\ref{partition}, $D\in\mathcal{C}(G)\cup\mathcal{S}(G)\cup\mathcal{P}(G)$. 
If $D\in\mathcal{P}(G)$, then $D=C$ because $D\in\mathcal{M}_{{\rm min}}(G)$ 
by Lemma~\ref{min}. 
So assume that $D\in\mathcal{C}(G)\cup\mathcal{S}(G)$. 
As $C$ is a clique in $G$, $D\in\mathcal{C}(G)$ and $C=D$ by the maximality of 
$C\in\mathcal{C}(G)$. 

Second, assume that $C\in\mathcal{P}(G)$. 
By Lemma~\ref{min}, $C\in\mathcal{M}_{{\rm min}}(G)$ and hence $c\approx_G d$ for any 
$c\neq d\in C$. 
So there is $D\in\mathfrak{M}(G)$ such that $D\supseteq C$. 
By Lemma~\ref{partition}, $D\in\mathcal{C}(G)\cup\mathcal{S}(G)\cup\mathcal{P}(G)$. 
As $G[C]$ is prime, $C$ is not included in a clique or a stable set in $G$. 
Therefore $D\in\mathcal{P}(G)$. 
By Lemma~\ref{min}, $D\in\mathcal{M}_{{\rm min}}(G)$ and hence $C=D$. 
\epr

Given a graph $G$, set $\mathcal{I}(G)=\mathfrak{M}(G)\setminus (\mathcal{C}(G)\cup\mathcal{S}(G)\cup\mathcal{P}(G))$ and $I(G)=\{v\in V(G):\{v\}\in\mathcal{I}(G)\}$. 
In Figure~\ref{fig1}, $I(G)=\{a,b\}$. 

\begin{rem}\label{base}
Given a graph $G$, consider $M\in\mathcal{M}(G)$. 
There exists $N\in\mathcal{M}_{{\rm min}}(G)$ such that $N\subseteq M$. 
By considering $p\neq q\in N$, we obtain that there exist $p\neq q\in M$ such that 
$p\approx_Gq$. By Corollary~\ref{cpartition}, there is 
$C\in\mathcal{C}(G)\cup\mathcal{S}(G)\cup\mathcal{P}(G)$ such that $|C\cap M|\geq 2$. 
\end{rem}

Let $G$ be a graph. 
If $\omega_M(G)\geq 2$, that is, $\mathcal{C}(G)\neq\emptyset$, then 
$\omega_M(G)=\max(\{\left|C\right|:C\in\mathcal{C}(G)\})$. 
Similarly 
$\alpha_M(G)=\max(\{\left|C\right|:C\in\mathcal{S}(G)\})$ if $\alpha_M(G)\geq 2$. 
Consequently, if $m(G)\geq 2$, that is, $\mathcal{C}(G)\cup\mathcal{S}(G)\neq\emptyset$, then 
$$m(G)=\max(\{\left|C\right|:C\in\mathcal{C}(G)\cup\mathcal{S}(G)\}).$$ 
Given a graph $G$, Sabidussi \cite{S59} introduced the following equivalence relation ${\rm Sab}_G$ on $V(G)$.  Given $u,v\in V(G)$, 
$u\hspace{1mm}{\rm Sab}_G\hspace{0.8mm}v$ if $N_G(u)=N_G(v)$. 
If $\alpha_M(G)\geq 2$, then $\mathcal{S}(G)$ is the family of the equivalence classes of ${\rm Sab}_G$ which are not singletons. 

We complete the section with another equivalence relation induced by the modular decomposition tree. 
It is used by Giakoumakis and Olariu \cite{GO07} to construct a minimal prime extension of a graph. 
Given a graph $G$, consider the equivalence relation $\leftrightarrow_G$ defined on $V(G)$ by: 
given $v,w\in V(G)$, $v\leftrightarrow_G w$ if 
$\{v\}\!\!\twoheaduparrow=\{w\}\!\!\twoheaduparrow$. 
The set of the equivalence classes of $\leftrightarrow_G$ is denoted by 
$\mathfrak{S}(G)$. 
In Figure~\ref{fig1}, $\mathfrak{S}(G)=\{C,S,\{a,b\}\}$. 

Precisely, to construct a minimal prime extension of a graph $G$, Giakoumakis and Olariu \cite{GO07} use only the elements of 
$\mathfrak{S}(G)\cap\mathcal{M}(G)$. 
The remainder of the section is mainly devoted to the study of $\mathfrak{S}(G)\cap\mathcal{M}(G)$ (see Proposition~\ref{summarize}). 

\begin{lem}\label{equivalence}
Let $G$ be a graph.
\begin{enumerate}
\item Let $M\in\mathbb{S}(G)$ with $|M|\geq 2$. 
If $\{N\in\Pi(G[M]):|N|=1\}\neq\emptyset$, then 
$\{m\in M:\{m\}\in\Pi(G[M])\}\in\mathfrak{S}(G)$. 
\item Let $C\in\mathfrak{S}(G)$ with $|C|\geq 2$. For every $c\in C$, $C\!\!\uparrow=\{c\}\!\!\twoheaduparrow$ and 
$C=\{c\in C\!\!\uparrow:\{c\}\in\Pi(G[C\!\!\uparrow])\}$. 
\item For each $M\in\mathbb{S}(G)$ with $|M|\geq 2$, there is $C\in\mathfrak{S}(G)\cap\mathbb{S}(G)$ such that $|C|\geq 2$ and $C\subseteq M$. 
\end{enumerate}
\end{lem}

\bpr
The first assertion follows from the definition of $\leftrightarrow_G$. 
For the second, consider $C\in\mathfrak{S}(G)$ with $|C|\geq 2$. 
For $c,d\in C$, we have $\{c\}\!\!\twoheaduparrow=\{d\}\!\!\twoheaduparrow$. 
Given $c_0\in C$, we obtain $C\subseteq\{c_0\}\!\!\twoheaduparrow$ and hence $C\!\!\uparrow\subseteq\{c_0\}\!\!\twoheaduparrow$. 
As $|C|\geq 2$, we have also $\{c_0\}\subsetneq C\!\!\uparrow$ and hence $\{c_0\}\!\!\twoheaduparrow\subseteq C\!\!\uparrow$. 
Thus $C\!\!\uparrow=\{c\}\!\!\twoheaduparrow$ for every $c\in C$. 
It follows from the definition of $\leftrightarrow_G$ that for each $v\in V(G)$, $v\in C$ if and only if $\{v\}\!\!\twoheaduparrow=C\!\!\uparrow$. 
Furthermore it follows from the second assertion of Proposition~\ref{strong} that for each $d\in C\!\!\uparrow$, 
$\{d\}\!\!\twoheaduparrow=C\!\!\uparrow$ if and only if $\{d\}\in\Pi(G[C\!\!\uparrow])$. 
Therefore $C=\{c\in C\!\!\uparrow:\{c\}\in\Pi(G[C\!\!\uparrow])\}$. 

For the third assertion, consider a strong module $M$ of $G$ with $|M|\geq 2$. 
Let $N$ be a minimal strong module of $G$ under inclusion among the strong modules $M'$ of $G$ such that $|M'|\geq 2$ and 
$M'\subseteq M$. 
By the minimality of $N$, we obtain $\Pi(G[N])=\{\{n\}:n\in N\}$. 
By the first assertion, $N\in\mathfrak{S}(G)$. 
\epr

\begin{prop}\label{compar}
For a graph $G$, 
\begin{enumerate}
\item $\mathcal{P}(G)=\{C\in\mathfrak{S}(G)\cap\mathbb{S}(G):\lambda_G(C)=\sqcup\}$,
\item $\mathcal{C}(G)=\{C\in\mathfrak{S}(G):|C|\geq 2,\lambda_G(C\!\!\uparrow)=\filledmedsquare\}$,
\item $\mathcal{S}(G)=\{C\in\mathfrak{S}(G):|C|\geq 2,\lambda_G(C\!\!\uparrow)=\medsquare\}$.
\end{enumerate}
\end{prop}

\bpr
For the first assertion, consider $P\in\mathcal{P}(G)$. 
By Lemma~\ref{prime}, $P\in\mathbb{S}(G)$. 
By the third assertion of Lemma~\ref{equivalence}, there is $C\in\mathfrak{S}(G)\cap\mathbb{S}(G)$ such that $|C|\geq 2$ 
and $C\subseteq P$. 
As $P\in\mathcal{M}_{{\rm min}}(G)$ by Lemma~\ref{min}, $C=P$. 
Since $G[P]$ is prime, $\Pi(G[P])=\{\{p\}:p\in P\}$. 
By the first assertion of Lemma~\ref{equivalence}, $P\in\mathfrak{S}(G)$. 
Moreover it follows from the first assertion of Proposition~\ref{pquotient} that $G[P]/\Pi(G[P])$ is prime. 
Thus $\lambda_G(P)=\sqcup$. 

Conversely, consider $C\in\mathfrak{S}(G)\cap\mathbb{S}(G)$ such that $\lambda_G(C)=\sqcup$. 
Clearly $C\!\!\uparrow=C$ because $C\in\mathbb{S}(G)$. 
As $C\in\mathfrak{S}(G)$, it follows from the second assertion of Lemma~\ref{equivalence} that $\Pi(G[C])=\{\{c\}:c\in C\}$. 
Furthermore we have $G[C]/\Pi(G[C])$ is prime because $\lambda_G(C)=\sqcup$. 
Thus $G[C]$ is prime as well by the first assertion of Proposition~\ref{pquotient}. 
Therefore $C\in\mathcal{P}(G)$. 

For the second assertion, consider $C\in\mathcal{C}(G)$. 
Denote by $Q$ the family of $M\in\Pi(G[C\!\!\uparrow])$ such that $M\cap C\neq\emptyset$. 
For each $M\in Q$, $C\setminus M\neq\emptyset$ because $M\subsetneq C\!\!\uparrow$. 
As $M$ is a strong module of $G$ by Proposition~\ref{strong}, we obtain $M\subsetneq C$. 
Thus $|Q|\geq 2$ and $C=\cup Q$. 
Given $M\in Q$, consider $N\in Q\setminus\{M\}$. 
For $m\in M$ and $n\in N$, $\{m,n\}$ is a module of $G[C]$ because $G[C]$ is complete. 
By the second assertion of Proposition~\ref{properties}, $\{m,n\}$ is a module of $G$. 
Since $M$ is a strong module of $G$ such that $m\in M\cap\{m,n\}$ and $n\in\{m,n\}\setminus M$, we get $M\subseteq\{m,n\}$ and hence $M=\{m\}$. 
Thus $\{c\}\in\Pi(G[C\!\!\uparrow])$ for each $c\in C$. 
Set $D=\{c\in C\!\!\uparrow:\{c\}\in\Pi(G[C\!\!\uparrow])\}$. 
We have $C\subseteq D$ and $D\in\mathfrak{S}(G)$ by the first assertion of Lemma~\ref{equivalence}. 
Moreover, by the first assertion of Proposition~\ref{properties}, $C$ is a module of $G[C\!\!\uparrow]$. 
It follows from the second assertion of Proposition~\ref{pquotient} that $\{\{c\}:c\in C\}$ is a module of $G[C\!\!\uparrow]/\Pi(G[C\!\!\uparrow])$. 
Since $\{\{c\}:c\in C\}$ is a clique in $G[C\!\!\uparrow]/\Pi(G[C\!\!\uparrow])$, we obtain 
$\lambda_G(C\!\!\uparrow)=\filledmedsquare$. 
As $\lambda_G(C\!\!\uparrow)=\filledmedsquare$ and as $D=\{c\in C\!\!\uparrow:\{c\}\in\Pi(G[C\!\!\uparrow])\}$, 
$\{\{d\}:d\in D\}$ is a module of $G[C\!\!\uparrow]/\Pi(G[C\!\!\uparrow])$ and a clique in $G[C\!\!\uparrow]/\Pi(G[C\!\!\uparrow])$. 
Consequently $D$ is a clique in $G$ and $D$ is a module of $G[C\!\!\uparrow]$ by the last assertion of Proposition~\ref{pquotient}. 
By the second assertion of Proposition~\ref{properties}, $D$ is a module of $G$. 
Since $C\in\mathcal{C}(G)$, $C=D$. 

Conversely, consider $C\in\mathfrak{S}(G)$ such that $|C|\geq 2$ and $\lambda_G(C\!\!\uparrow)=\filledmedsquare$. 
By the second assertion of Lemma~\ref{equivalence}, $C=\{c\in C\!\!\uparrow:\{c\}\in\Pi(G[C\!\!\uparrow])\}$. 
Since $\lambda_G(C\!\!\uparrow)=\filledmedsquare$, $C$ is a clique in $G$ and, as above for $D$, it follows from 
Propositions~\ref{properties} and \ref{pquotient} that $C$ is a module of $G$. 
Thus there is $D\in\mathcal{C}(G)$ such that $C\subseteq D$. 
As already proved, $D\in\mathfrak{S}(G)$ and hence $C=D$. 

The third assertion follows from the second by interchanging $G$ and $\overline{G}$. 
\epr

It follows from Proposition~\ref{compar} that 

\begin{cor}
For a graph $G$,  
\begin{equation*}
\begin{cases}
\omega_\mathcal{M}(G)=\max(\{|C|:C\in\mathfrak{S}(G),|C|\geq 2,\lambda_G(C\!\!\uparrow)=\filledmedsquare\})\ \text{if}\ \omega_\mathcal{M}(G)\geq 2,\\
\alpha_\mathcal{M}(G)=
\max(\{|C|:C\in\mathfrak{S}(G),|C|\geq 2,\lambda_G(C\!\!\uparrow)=\medsquare\})\ \text{if}\ 
\alpha_\mathcal{M}(G)\geq 2.
\end{cases}
\end{equation*}
\end{cor}

The next is a simple consequence of Corollary~\ref{cpartition} and Proposition~\ref{compar}. 

\begin{cor}\label{ccompar}
Let $G$ be a graph. 
\begin{enumerate}
\item $\{C\in\mathfrak{M}(G):|C|\geq 2\}\subseteq\{C\in\mathfrak{S}(G):|C|\geq 2\}$. 
\item Given $C \in\mathfrak{S}(G)$ with $|C|\geq 2$, 
$$C\not\in\mathfrak{M}(G)\ \text{if and only if}\ C\neq C\!\!\uparrow\text{and}\ \lambda_G(C\!\!\uparrow)=\sqcup.$$
\item Let $C\in\mathfrak{S}(G)\setminus\mathfrak{M}(G)$ with $|C|\geq 2$. For each $c\in C$, 
$\{c\}\in\mathfrak{M}(G)\setminus\mathfrak{S}(G)$. 
\item $I(G)=\{v\in V(G):\{v\}\in\mathfrak{S}(G)\}
\cup (\bigcup_{C\in\mathfrak{S}(G)\setminus\mathfrak{M}(G)}\{\{c\}:c\in C\})$.
\end{enumerate}
\end{cor}

\begin{rem}
The second assertion of Corollary~\ref{ccompar} easily provides graphs $G$ such that 
$\mathfrak{S}(G)\setminus\mathfrak{M}(G)\neq\emptyset$. 
For instance, $\{a,b\}\in\mathfrak{S}(G)\setminus\mathfrak{M}(G)$ in Figure~\ref{fig1}.

Consider a prime graph $G_0$. 
Given $\alpha\not\in V(G_0)$ and $u\in V(G_0)$, consider a 1-extension $G$ of $G_0$ to $V(G_0)\cup\{\alpha\}$ such that 
$\{u,\alpha\}$ is a module of $G$, that is, $\alpha\in X_G(u)$ where $X=V(G_0)$. 
Clearly $\{u,\alpha\}$ is the single non-trivial module of $G$. 
Consequently
\begin{equation*}
\begin{cases}
\mathfrak{M}(G)=\{\{u,\alpha\}\}\cup\{\{v\}:v\in V(G_0)\setminus\{u\}\},\\
\mathbb{S}(G)=\{\{v\}:v\in V(G)\}\cup\{\{u,\alpha\},V(G)\}.\\
\end{cases}
\end{equation*}
Since $\mathbb{S}(G)=\{\{v\}:v\in V(G)\}\cup\{\{u,\alpha\},V(G)\}$, we get 
$$\mathfrak{S}(G)=\{\{u,\alpha\},V(G_0)\setminus\{u\}\}.$$
Thus $V(G_0)\setminus\{u\}\in\mathfrak{S}(G)\setminus\mathfrak{M}(G)$. 
\end{rem}

The following summarizes our comparison between $\mathfrak{M}(G)$ and $\mathfrak{S}(G)$ for a graph $G$. 

\begin{prop}\label{summarize}
For a graph $G$, 
\begin{align*}
\mathfrak{S}(G)\cap\mathcal{M}(G)&=\{C\in\mathfrak{M}(G):|C|\geq 2\}\\
&=\mathcal{C}(G)\cup\mathcal{S}(G)\cup\mathcal{P}(G).
\end{align*}
\end{prop}

\bpr
By Corollary~\ref{cpartition}, $\{C\in\mathfrak{M}(G):|C|\geq 2\}=\mathcal{C}(G)\cup\mathcal{S}(G)\cup\mathcal{P}(G)$. 
Furthermore, it follows from Proposition~\ref{compar} that 
$\mathcal{C}(G)\cup\mathcal{S}(G)\cup\mathcal{P}(G)\subseteq\mathfrak{S}(G)\cap\mathcal{M}(G)$. 
So consider $C\in\mathfrak{S}(G)\cap\mathcal{M}(G)$. 
If $\lambda_G(C\!\!\uparrow)=\filledmedsquare$ or $\medsquare$, then $C\in\mathcal{C}(G)\cup\mathcal{S}(G)$ by the last two assertions of Proposition~\ref{compar}. 
Thus assume that $\lambda_G(C\!\!\uparrow)=\sqcup$, that is, $G[C\!\!\uparrow]/\Pi(G[C\!\!\uparrow])$ is prime. 
By the second assertion of Lemma~\ref{equivalence}, $C=\{c\in C\!\!\uparrow:\{c\}\in\Pi(G[C\!\!\uparrow])\}$. 
Since $C$ is a module of $G$, $C$ is a module of $G[C\!\!\uparrow]$ by the first assertion of Proposition~\ref{properties}. 
Therefore $\{\{c\}:c\in C\}$ is a module of $G[C\!\!\uparrow]/\Pi(G[C\!\!\uparrow])$ by the second assertion of Proposition~\ref{pquotient}. 
As $G[C\!\!\uparrow]/\Pi(G[C\!\!\uparrow])$ is prime, $\{\{c\}:c\in C\}=\Pi(G[C\!\!\uparrow])$. 
Consequently, $C=C\!\!\uparrow$ and hence $C\in\mathcal{P}(G)$ by the first assertion of Proposition~\ref{compar}. 
\epr

Given a non-primitive and connected graph $G$, Giakoumakis and Olariu \cite[Theorem 3.9]{GO07} construct a minimal prime extension of 
$G$ by adding $|C|-1$ vertices for each $C\in\mathcal{C}(G)\cup\mathcal{S}(G)$ and one vertex for each element of $\mathcal{P}(G)$.

\section{Some prime extensions}

We use the next corollary to prove Theorem~\ref{mainone}. 

\begin{lem}\label{2Xstable}
Let $S$ and $S'$ be disjoint sets such that $|S'|=\lceil\log_2(|S|+1)\rceil\geq 2$, that is, $1\leq 2^{|S'|-1}\leq |S|< 2^{|S'|}$. 
Consider any graph $G$ defined on $V(G)=S\cup S'$ such that $S$ and $S'$ are stable sets in $G$. 
\begin{enumerate}
\item Assume that $2^{|S'|-1}< |S|< 2^{|S'|}$. Then, $G$ is prime if and only if 
\begin{equation}\label{N_G}
\text{the function $(N_G)_{\restriction S}$ is an injection from $S$ into $2^{S'}\setminus\{\emptyset\}$;}
\end{equation}
\item Assume that $2^{|S'|-1}=|S|$. Then, $G$ is prime if and only if \eqref{N_G} holds and for any  $s\in S$ and $s'\in S'$, 
\begin{equation*}
\text{if $d_G(s)=d_G(s')=1$, then $(s,s')_G=0$.}
\end{equation*}
\end{enumerate}
\end{lem}

\bpr
First, if there is $s\in S$ such that $N_G(s)=\emptyset$, then $s\in{\rm Iso}(G)$ and hence $V(G)\setminus\{s\}$ is a non-trivial module of $G$. 
Second, if there are $s\neq t\in S$ such that $N_G(s)=N_G(t)$, then 
$\{s,t\}$ is a non-trivial module of $G$. 
Third, if there are $s\in S$ and $s'\in S'$ such that $N_G(s)=\{s'\}$ and $N_G(s')=\{s\}$, then 
$\{s,s'\}$ is a non-trivial module of $G$. 

Conversely, assume that \eqref{N_G} holds. 
Consider a module $M$ of $G$ such that $\left|M\right|\geq 2$. 
We have to show that $M=V(G)$. 
As $(N_G)_{\restriction S}$ is injective, $M\not\subseteq S$, that is, $M\cap S'\neq\emptyset$. 

For a first contradiction, suppose that $M\subseteq S'$. 
Recall that for each $s\in S$, either $M\cap N_G(s)=\emptyset$ or $M\subseteq N_G(s)$. 
Given $m\in M$, 
consider the function $f:S\longrightarrow 2^{((S'\setminus M)\cup\{m\})}\setminus\{\emptyset\}$ defined by 
\begin{equation*}
f(s)=
\begin{cases}
\text{$N_G(s)$ if $M\cap N_G(s)=\emptyset$,}\\
\text{$N_G(s)\cup\{m\}$ if $M\subseteq N_G(s)$,}
\end{cases}
\end{equation*}
for every $s\in S$. 
Since $(N_G)_{\restriction S}$ is injective, $f$ is also and we would obtain that 
$\left|S\right|<2^{\left|S'\right|-1}$. 
It follows that $M\cap S\neq\emptyset$. 

For a second contradiction, suppose that $S'\setminus M\neq\emptyset$. 
We have $(S\cap M,S'\setminus M)_G=(S'\cap M,S'\setminus M)_G=0$. 
If $S\subseteq M$, then $(S,S'\setminus M)_G=0$ so that 
$(N_G)_{\restriction S}$ would be an injection from $S$ into $2^{S'\cap M}\setminus\{\emptyset\}$ which contradicts $|S|\geq 2^{|S'|-1}$. 
Thus $S\setminus M\neq\emptyset$. 
We obtain $(S\setminus M,S'\cap M)_G=(S\setminus M,S\cap M)_G=0$. 
As $(S\cap M,S'\setminus M)_G=0$, 
$(N_G)_{\restriction (S\cap M)}:S\cap M\longrightarrow 2^{S'\cap M}\setminus\{\emptyset\}$ and hence $\left|S\cap M\right|\leq 2^{\left|S'\cap M\right|}-1$. 
Since $(S\setminus M,S'\cap M)_G=0$, $(N_G)_{\restriction (S\setminus M)}:S\setminus M\longrightarrow 2^{S'\setminus M}\setminus\{\emptyset\}$ and hence $\left|S\setminus M\right|\leq 2^{\left|S'\setminus M\right|}-1$. 
Therefore $\left|S\right|\leq 2^{\left|S'\cap M\right|}+2^{\left|S'\setminus M\right|}-2\leq 2^{\left|S'\right|-1}$. 
As $\left|S\right|\geq 2^{\left|S'\right|-1}$, we obtain 
$\left|S\right|=2^{\left|S'\right|-1}$ and $2^{\left|S'\cap M\right|}+2^{\left|S'\setminus M\right|}=2+2^{\left|S'\right|-1}$ 
so that $\min(\left|S'\cap M\right|,\left|S'\setminus M\right|)=1$. 
For instance, assume that $\left|S'\cap M\right|=1$. 
Since $\left|S\cap M\right|\leq 2^{\left|S'\cap M\right|}-1$, $\left|S\cap M\right|=1$. 
There exist $s\in S$ and $s'\in S'$ such that $M=\{s,s'\}$. 
By what precedes, $(s,S'\setminus\{s'\})_G=(s',S\setminus\{s\})_G=0$. 
As $N_G(s)\neq\emptyset$, we would obtain $N_G(s)=\{s'\}$, 
$N_G(s')=\{s\}$ and $(s,s')_G=1$. 
It follows that $S'\subseteq M$. 

Lastly, suppose that $S\setminus M\neq\emptyset$. 
For each $s\in S\setminus M\neq\emptyset$, we would have $(s,S')_G=(s,S\cap M)_G=0$ and hence 
$N_G(s)=\emptyset$. 
It follows that $S\subseteq M$ and $M=S\cup S'$. 
\epr

\begin{cor}\label{c2Xstable}
Let $S$ and $S'$ be disjoint sets such that $|S|\geq 3$ and $|S'|=\lceil\log_2(|S|+1)\rceil$. 
There exists a prime graph $G$ defined on $V(G)=S\cup S'$ and satisfying  
\begin{itemize}
\item $S$ and $S'$ are stable sets in $G$;
\item $(N_G)_{\restriction S}:S\longrightarrow 2^{S'}\setminus\{\emptyset\}$ is injective;
\item there exists an injection $\varphi_{S'}:S'\longrightarrow S$ such that 
$N_G(\varphi_{S'}(s'))=S'\setminus\{s'\}$ for each $s'\in S'$.
\end{itemize}
\end{cor}

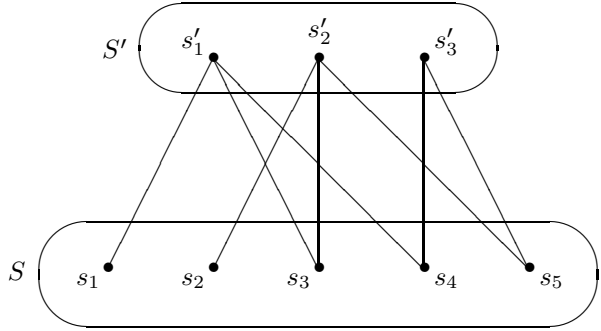
\begin{figure}[h]
\begin{center}
\setlength{\unitlength}{0.7cm}
\begin{picture}(14,7)

\put(3,2){$\bullet$}\put(2.5,1.8){$s_1$}
\put(5,2){$\bullet$}\put(4.5,1.8){$s_2$}
\put(7,2){$\bullet$}\put(6.5,1.8){$s_3$}
\put(9,2){$\bullet$}\put(9.3,1.8){$s_4$}
\put(11,2){$\bullet$}\put(11.3,1.8){$s_5$}

\put(7.1,2){\oval(10.6,2)}
\put(1.2,1.9){$S$}

\put(5,6){$\bullet$}\put(4.5,6.3){$s'_1$}
\put(7,6){$\bullet$}\put(6.9,6.5){$s'_2$}
\put(9,6){$\bullet$}\put(9.3,6.3){$s'_3$}

\put(7.1,6.3){\oval(6.8,1.7)}
\put(3,6.1){$S'$}

\put(3.1,2.1){\line(1,2){2}}
\put(5.1,2.1){\line(1,2){2}}
\put(7.1,2.1){\line(0,1){4}}
\put(7.1,2.1){\line(-1,2){2}}
\put(9.1,2.1){\line(0,1){4}}
\put(9.1,2.1){\line(-1,1){4}}
\put(11.1,2.1){\line(-1,1){4}}
\put(11.1,2.1){\line(-1,2){2}}

\end{picture}
\caption{Corollary~\ref{c2Xstable} when $S=\{s_1,\ldots,s_5\}$, $S'=\{s'_1,s'_2,s'_3\}$ 
with $\varphi_{S'}(s'_1)=s_5$, $\varphi_{S'}(s'_2)=s_4$ and $\varphi_{S'}(s'_3)=s_3$.\label{fig2}}
\end{center}
\end{figure}

\bpr 
(See Figure~\ref{fig2}.) 
As $|S'|=\lceil\log_2(|S|+1)\rceil$, we have $2^{|S'|-1}\leq |S|< 2^{|S'|}$ and hence 
$|S'|\leq|S|$. 
Thus there exists a bijection $\psi_{S'}$ from $S'$ onto $S''\subseteq S$. 
Consider the injection $f_{S''}:S''\longrightarrow 2^{S'}\setminus\{\emptyset\}$ defined by 
$s''\mapsto S'\setminus\{(\psi_{S'})^{-1}(s'')\}$. 
Let $f_S$ be any injection from $S$ into $2^{S'}\setminus\{\emptyset\}$ such that 
$(f_S)_{\restriction S''}=f_{S''}$. 
Lastly, consider the graph $G$ defined on $V(G)=S\cup S'$ such that $S$ and $S'$ are stable sets  in $G$ and $(N_G)_{\restriction S}=f_S$. 
Before applying Lemma~\ref{2Xstable}, assume that $|S|=2^{|S'|-1}$. 
Since $|S|\geq 3$, $|S'|\geq 3$. 
For each $s'\in S'$, there are $t'\neq u'\in S'\setminus\{s'\}$. 
We obtain $N_G(\psi_{S'}(t'))=S'\setminus\{t'\}$ and $N_G(\psi_{S'}(u'))=S'\setminus\{u'\}$. 
Therefore $s'\in N_G(\psi_{S'}(t'))\cap N_G(\psi_{S'}(u'))$. 
It follows that $d_G(s')\geq 2$ for every $s'\in S'$. 
By Lemma~\ref{2Xstable}, $G$ is prime. 
\epr

We use the following two results to prove Theorem~\ref{mainone} when $\mathcal{P}(G)\neq\emptyset$. 
Given a graph $G$, consider $X\subsetneq V(G)$ such that $G[X]$ is prime. 
We utilize the following subsets of $V(G)\setminus X$ 
(for instance, see~\cite[Lemma~5.1]{ER90})
\begin{itemize}
\item ${\rm Ext}_G(X)$ is the set of $v\in V(G)\setminus X$ such that $G[X\cup\{v\}]$ is prime;
\item $\langle X\rangle_G$ is the set of $v\in V(G)\setminus X$ such that $X$ is a module of $G[X\cup\{v\}]$;
\item for $u\in X$, $X_G(u)$ is the set of $v\in V(G)\setminus X$ such that $\{u,x\}$ is a module of $G[X\cup\{v\}]$. 
\end{itemize}

The family $\{{\rm Ext}_G(X),\langle X\rangle_G\}\cup\{X_G(u):u\in X\}$ is denoted by $p_{G[X]}$. 
The next lemma follows from Proposition~\ref{properties}. 

\begin{lem}\label{G[X]}
Given a graph $G$, consider $X\subsetneq V(G)$ such that $G[X]$ is prime. 
The family $p_{G[X]}$ is a partition of $V(G)\setminus X$. 
Moreover, 
for each module $M$ of $G$, one and only one of the following holds
\begin{itemize}
\item $X\subseteq M$ and $V(G)\setminus M\subseteq\langle X\rangle_G$; 
\item there is a unique $u\in X$ such that $M\cap X=\{u\}$ and $M\setminus\{u\}\subseteq X_G(u)$;
\item $M\cap X=\emptyset$ and $M$ is included in an element of $p_{G[X]}$. Moreover, for $v,w\in M$, the function 
$X\cup\{v\}\longrightarrow X\cup\{w\}$, defined by $v\mapsto w$ and $u\mapsto u$ for $u\in X$, is an isomorphism from $G[X\cup\{v\}]$ onto 
$G[X\cup\{w\}]$. 
\end{itemize}
\end{lem}

\begin{lem}\label{1extension}
Let $G$ be a prime graph. For every $\alpha\not\in V(G)$, there are $$2^{|V(G)|}-2|V(G)|-2$$ distinct prime 1-extensions of $G$ to 
$V(G)\cup\{\alpha\}$. 
\end{lem}

\bpr
Consider any graph $H$ defined on $V(H)=V(G)\cup\{\alpha\}$ such that $H[V(G)]=G$ and 
$$N_H(\alpha)\in 2^{V(G)}\setminus (\{\emptyset,V(G)\}\cup\{N_G(v):v\in V(G)\}\cup\{N_G(v)\cup\{v\}:v\in V(G)\}).$$
We verify that $H$ is prime. 
Set $X=V(G)$. We have $H[X]=G$ is prime. 
If $\alpha\in\langle X\rangle_H$, then $N_H(\alpha)=\emptyset$ or $V(G)$. 
Thus $\alpha\not\in\langle X\rangle_H$. 
If there is $v\in V(G)$ such that $\alpha\in X_H(v)$, then $N_H(\alpha)=N_G(v)$ or $N_G(v)\cup\{v\}$. 
Therefore $\alpha\not\in X_H(v)$ for every $v\in V(G)$. 
It follows from Lemma~\ref{G[X]} that $\alpha\in{\rm Ext}_H(X)$, that is, $H$ is prime. 
Consequently the number of prime 1-extensions of $G$ to $V(G)\cup\{\alpha\}$ equals 
$$|2^{V(G)}\setminus (\{\emptyset,V(G)\}\cup\{N_G(v):v\in V(G)\}\cup\{N_G(v)\cup\{v\}:v\in V(G)\})|.$$

Clearly $\emptyset\not\in\{N_G(v)\cup\{v\}:v\in V(G)\}$, $V(G)\not\in\{N_G(v):v\in V(G)\}$ and 
$\{N_G(v):v\in V(G)\}\cap\{N_G(v)\cup\{v\}:v\in V(G)\}=\emptyset$. 
Moreover, 
if there is $v\in V(G)$ such that $N_G(v)=\emptyset$ or $V(G)\setminus\{v\}$, then $V(G)\setminus\{v\}$ would be a non-trivial module of $G$. 
If there are $v\neq w\in V(G)$ such that $N_G(v)=N_G(w)$ or 
$N_G(v)\cup\{v\}=N_G(w)\cup\{w\}$, then $\{v,w\}$ would be a non-trivial module of $G$. 
As $G$ is prime, $V(G)\not\in\{N_G(v)\cup\{v\}:v\in V(G)\}$, $\emptyset\not\in\{N_G(v):v\in V(G)\}$ and for $v\neq w\in V(G)$, we have 
$N_G(v)\neq N_G(w)$ and $N_G(v)\cup\{v\}\neq N_G(w)\cup\{w\}$. 
Therefore $|2^{V(G)}\setminus (\{\emptyset,V(G)\}\cup\{N_G(v):v\in V(G)\}\cup\{N_G(v)\cup\{v\}:v\in V(G)\})|=2^{|V(G)|}-2|V(G)|-2$. 
\epr

\section{Proofs of Theorems~\ref{mainone} and \ref{maintwo}}

\begin{figure}[h]
\begin{center}
\setlength{\unitlength}{0.7cm}
\begin{picture}(14,8)

\put(0,3){$\bullet$}\put(-0.5,2.8){$s_1$}
\put(2,3){$\bullet$}\put(1.5,2.8){$s_2$}
\put(4,3){$\bullet$}\put(3.5,2.8){$s_3$}
\put(6,3){$\bullet$}\put(6.3,2.8){$s_4$}
\put(8,3){$\bullet$}\put(8.3,2.8){$s_5$}

\put(0.1,3.1){\line(3,2){6}}
\put(2.1,3.1){\line(3,2){6}}
\put(4.1,3.1){\line(1,2){2}}
\put(4.1,3.1){\line(1,1){4}}
\put(6.1,3.1){\line(0,1){4}}
\put(6.1,3.1){\line(1,1){4}}
\put(8.1,3.1){\line(0,1){4}}
\put(8.1,3.1){\line(1,2){2}}

\put(4.1,3){\oval(9.6,2)}
\put(0.1,2.2){$S$}

\put(4.1,0.1){\line(-4,3){4}}
\put(4.1,0.1){\line(-2,3){2}}
\put(4.1,0.1){\line(0,1){3}}
\put(4.1,0.1){\line(2,3){2}}
\put(4.1,0.1){\line(4,3){4}}

\put(12,3){$\bullet$}\put(11.6,2.8){$c_1$}
\put(14,3){$\bullet$}\put(14.3,2.8){$c_2$}
\put(13,1.5){$\bullet$}\put(12.87,2.1){$c_3$}

\put(13.1,2.5){\oval(3.4,2.2)}
\put(14.1,1.7){$C$}

\put(13.1,0.1){\line(0,1){1.5}}
\put(13.1,0.1){\line(-1,3){1}}
\put(13.1,0.1){\line(1,3){1}}

\put(13.1,0.1){\line(-1,1){7}}

\put(13.1,1.6){\line(-2,3){1}}
\put(13.1,1.6){\line(2,3){1}}
\put(12.1,3.1){\line(1,0){2}}

\put(12.1,3.1){\line(-1,1){4}}
\put(14.1,3.1){\line(-1,1){4}}

\put(4,0){$\bullet$}\put(3.6,-0.2){$a$}
\put(13,0){$\bullet$}\put(13.4,-0.2){$b$}
\put(4.1,0.1){\line(1,0){9}}

\put(7.2,2){\oval(16,4.7)}
\put(0,0){$V(G)$}

\put(6,7){$\bullet$}\put(5.9,7.5){$s'_1$}
\put(8,7){$\bullet$}\put(7.9,7.5){$s'_2$}
\put(10,7){$\bullet$}\put(9.9,7.5){$s'_3$}

\put(8.1,7.3){\oval(5,2)}
\put(4.7,7.1){$S'$}

\end{picture}
\caption{Theorem~\ref{mainone} for the graph depicted in Figure~\ref{fig1}.\label{fig3}}
\end{center}
\end{figure}
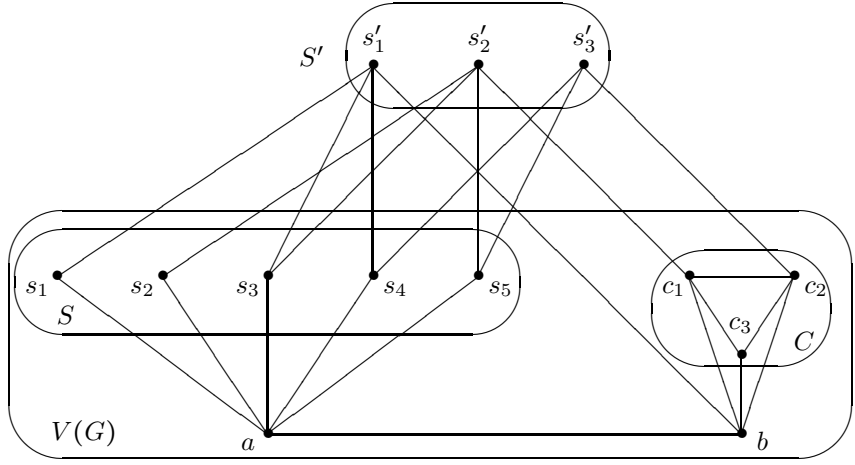

\bprthone
Consider a graph $G$ such that $|m(G)|\geq 2$. 
To begin, we show that $p(G)\geq\lceil\log_2(m(G))\rceil$. 
By interchanging $G$ and $\overline{G}$, assume that there exists $S\in\mathcal{S}(G)$ with 
$\left|S\right|=m(G)$. 
Given an integer $p<\log_2(m(G))$, 
consider any $p$-extension $H$ of $G$. 
We must prove that $H$ is not prime. 
We have $2^{\left|V(H)\setminus V(G)\right|}<\left|S\right|$ so that the function 
$S\longrightarrow 2^{V(H)\setminus V(G)}$, defined by $s\mapsto N_H(s)\cap (V(H)\setminus V(G))$, is not injective. 
Thus there are $s\neq t\in S$ such that 
$N_H(s)\cap (V(H)\setminus V(G))=N_H(t)\cap (V(H)\setminus V(G))$. 
In other words, $v\sim_H\{s,t\}$ for every $v\in V(H)\setminus V(G)$. 
As $S$ is a module of $G$, we have $v\sim_G\{s,t\}$, that is, $v\sim_H\{s,t\}$ 
for every $v\in V(G)\setminus S$. 
Since $S$ is a stable set in $G$, $(v,s)_H=(v,s)_G=(v,t)_G=(v,t)_H$. 
Therefore $\{s,t\}$ is a module of $H$ and $H$ is not prime.

To \ prove\  that\  $p(G)\leq\lceil\log_2(m(G)+1)\rceil$, we\  must\  construct\  a\  prime $\lceil\log_2($ $m(G)+1)\rceil$-extension $H$ of $G$. 
Let $S'$ be a set such that $S'\cap V(G)=\emptyset$ and 
$\left|S'\right|=\lceil\log_2(m(G)+1)\rceil$. 
Let $s'_1\in S'$. 
Consider $S_0\in\mathcal{C}(G)\cup\mathcal{S}(G)$ such that $\left|S_0\right|=m(G)$, that is, 
$2^{\left|S'\right|-1}\leq \left|S_0\right|<2^{\left|S'\right|}$. 
By interchanging $G$ and $\overline{G}$, we can assume that $S_0$ is a clique or a stable set in $G$. 
(In Figure~\ref{fig3}, $m(G)=5$, $S_0=S=\{s_1,\ldots,s_5\}\in\mathcal{S}(G)$ and 
$S'=\{s'_1,s'_2,s'_3\}$.)

We consider any graph $H$ defined on $V(G)\cup S'$ and satisfying the following. 
\begin{enumerate}
\item $S'$ is a stable set in $H$.
\item The subgraph $H[S_0\cup S']$ of $H$ is defined as follows
\begin{itemize}
\item Assume that $\left|S_0\right|=2$. We require that the subgraph $H[S_0\cup S']$ of $H$ is a path on 4 vertices and $S_0$ is a clique in $H[S_0\cup S']$. Thus 
\begin{equation}\label{S0=2}
d_{H[S_0\cup S']}(s')=1\ \text{for each $s'\in S'$.} 
\end{equation}
\item Assume that $\left|S_0\right|\geq 3$. By Corollary~\ref{c2Xstable}, we can consider for 
$H[S_0\cup S']$ a prime graph defined on $S_0\cup S'$ such that 
\begin{equation}\label{S0>2}
\begin{cases}
\text{$S_0$ is a stable set in $H[S_0\cup S']$;}\\
\text{$(N_{H[S_0\cup S']})_{\restriction S_0}:S_0\longrightarrow 2^{S'}\setminus\{\emptyset\}$ is injective;}\\
\text{there exists an injection $\varphi_{S'}:S'\longrightarrow S_0$ with}\\
\text{\hspace{15mm}$N_{H[S_0\cup S']}(\varphi_{S'}($ $s'))=S'\setminus\{s'\}$ for each $s'\in S'$.}
\end{cases}
\end{equation}
(In Figure~\ref{fig3}, the subgraph $H[S_0\cup S']$ is also depicted in Figure~\ref{fig2}.)
\end{itemize}
Set $X=S_0\cup S'$. In both cases, $H[X]$ is prime. 
\item Let $C\in\mathcal{C}(G)\setminus\{S_0\}$. We have $\left|C\right|<2^{\left|S'\right|}$. We consider for $H[C\cup S']$ a graph such that $C$ is a clique in $H[C\cup S']$ and 
\begin{equation}\label{C}
\left.\begin{aligned}
\hspace{-5mm}f_C:&\ C\longrightarrow 2^{S'}\\
&\ \ c\longmapsto N_{H[C\cup S']}(c)\cap S'
\end{aligned}\right\}\ \text{is an injection from $C$ into $2^{S'}\setminus\{S'\}$.}
\end{equation}
(In Figure~\ref{fig3}, we have $C=\{c_1,c_2,c_3\}\in\mathcal{C}(G)$ and $f_C$ is defined by 
$$c_1\mapsto\{s'_2\},\ c_2\mapsto\{s'_3\},\ c_3\mapsto\emptyset.)$$
\item Let $S\in\mathcal{S}(G)\setminus\{S_0\}$. We have $\left|S\right|<2^{\left|S'\right|}$. We consider for $H[S\cup S']$ a graph such that $S$ is a stable set in $H[S\cup S']$ and 
\begin{equation}\label{S}
\left.\begin{aligned}
f_S:&\ S\longrightarrow 2^{S'}\\
&\ \ s\longmapsto N_{H[S\cup S']}(s)
\end{aligned}\right\}\ \text{is an injection from $S$ into $2^{S'}\setminus\{\emptyset\}$.}
\end{equation}
\item Let $P\in\mathcal{P}(G)$. As $G[P]$ is prime, it follows from Lemma~\ref{1extension} that  $G[P]$ admits a prime 1-extension. We consider for $H[P\cup S']$ a graph such that 
\begin{equation}\label{P}
\begin{cases}
H[P]=G[P]\ \text{is prime,}\\
{\rm Ext}_{H[P\cup S']}(P)=S'\ \text{by using Lemma~\ref{1extension}.} 
\end{cases}
\end{equation}
\item Let $v\in I(G)$. Since $S_0$ is a module of $G$ such that $v\not\in S_0$, there is $i\in\{0,1\}$ such that $(v,S_0)_G=i$. We consider for $H[\{v\}\cup S']$ a graph such that 
\begin{equation}\label{I}
(v,s'_1)_{H[\{v\}\cup S']}\neq i.
\end{equation}
(In Figure~\ref{fig3}, we have $S_0=S$ and $I(G)=\{a,b\}$ with
\begin{equation*}
\begin{cases}
\text{$(a,S)_G=1$ and $(a,s'_1)_{H[\{a\}\cup S']}=0$,}\\
\text{$(b,S)_G=0$ and $(b,s'_1)_{H[\{b\}\cup S']}=1$.)}
\end{cases}
\end{equation*}
\end{enumerate}
In the construction above, we have $H[N]=G[N]$ for each $N\in\mathfrak{M}(G)$. Thus we can also assume that $H[V(G)]=G$. 

We begin with the following observation. 
For every module $M$ of $H$ such that $|M|\geq 2$, we have $M\cap S'\neq\emptyset$. 
Otherwise suppose that $M\subseteq V(G)$. 
By the first assertion of Proposition~\ref{properties}, $M\in\mathcal{M}(G)$. 
By Remark~\ref{base}, there is $C\in\mathcal{C}(G)\cup\mathcal{S}(G)\cup\mathcal{P}(G)$ such that $|C\cap M|\geq 2$. 
As $M$ is a module of $H$, it follows from the first assertion of Proposition~\ref{properties} that $C\cap M$ is a module of $H[C\cup S']$. 
Since $H[S_0\cup S']$ is prime, $C\neq S_0$. 
If $C\in (\mathcal{C}(G)\cup\mathcal{S}(G))\setminus\{S_0\}$, then 
$(f_C)_{\restriction (C\cap M)}$ would be constant which contradicts \eqref{C} and \eqref{S}. 
Lastly, suppose that $C\in\mathcal{P}(G)$. 
As $C\cap M$ is a module of $H[C]=G[C]$ by the first assertion of 
Proposition~\ref{properties}, it follows from Lemma~\ref{prime} that $C\cap M=C$. 
Thus $C$ would be a module of $H[C\cup S']$ 
so that 
$\langle C\rangle_{H[C\cup S']}=S'$ 
which contradicts \eqref{P} and Lemma~\ref{G[X]}. 
Consequently, $M\cap S'\neq\emptyset$ for every module $M$ of $H$ such that $|M|\geq 2$. 

Now, we prove that $H$ is prime. 
Consider a module $M$ of $H$ such that $|M|\geq 2$. 
We have to show that $M=V(H)$. 
As observed above, $M\cap S'\neq\emptyset$. 
By Lemma~\ref{G[X]}, either there is $s'\in S'$ such that $M\cap X=\{s'\}$ or $X\subseteq M$. 

For a contradiction, suppose that there is $s'\in S'$ such that $M\cap X=\{s'\}$. 
By the first assertion of Proposition~\ref{properties}, $M\setminus\{s'\}$ is a module of $G$. 
By the last assertion of Proposition~\ref{properties}, there is $i\in\{0,1\}$ such that $(M\setminus\{s'\},S_0)_G=i$. 
Thus $(s',S_0)_{H[S_0\cup S']}=i$. 
Since $S'$ is a stable set in $H[S_0\cup S']$, we have $d_{H[S_0\cup S']}(s')=0$ or $d_{H[S_0\cup S']}(s')\geq 2$ so that \eqref{S0=2} does not hold. 
Therefore $|S_0|\geq 3$. 
Let $t'\in S'\setminus\{s'\}$. 
By \eqref{S0>2}, $(\varphi_{S'}(s'),s')_{H[S_0\cup S']}=0$ and 
$(\varphi_{S'}(t'),s')_{H[S_0\cup S']}=1$ which contradicts $(s',S_0)_{H[S_0\cup S']}=i$. 
It follows that $X\subseteq M$. 

First, consider $C\in\mathcal{C}(G)\setminus\{S_0\}$. 
Suppose for a contradiction that $C\cap M=\emptyset$. 
Thus $C\cap S_0=\emptyset$ and it follows from the last assertion of Proposition~\ref{properties} that there is $i\in\{0,1\}$ such that $(C,S_0)_G=i$. 
As $C\cap M=\emptyset$, we obtain $(C,M)_H=i$. 
In particular $(C,S')_{H[C\cup S']}=i$ and $f_C$ would be constant which contradicts \eqref{C}. 
Therefore $C\cap M\neq\emptyset$. 
Suppose for a contradiction that $C\setminus M\neq\emptyset$ and consider $c\in C\setminus M$. 
Since $C$ is a clique of $G$, $(c,C\cap M)_G=1$. 
Hence $(c,M)_H=1$ and in particular $(c,S')_{H[C\cup S']}=1$ which contradicts $f_C(c)\neq S'$. 
It follows that $C\subseteq M$. 
Similarly $S\subseteq M$ for every $S\in\mathcal{S}(G)\setminus\{S_0\}$. 

Second, consider $P\in\mathcal{P}(G)$. 
By the first assertion of Proposition~\ref{properties}, $M\cap V(G)$ is a module of $G$. 
As $M\cap V(G)\supseteq S_0$, it follows from Lemma~\ref{prime} that either 
$(M\cap V(G))\cap P=\emptyset$ or $P\subseteq M\cap V(G)$. 
Suppose for a contradiction that $(M\cap V(G))\cap P=\emptyset$. 
By the last assertion of Proposition~\ref{properties}, there is $i\in\{0,1\}$ such that $(P,S_0)_G=i$. 
As $S_0\subseteq M$ and $M\cap P=\emptyset$, we obtain $(P,M)_H=i$. 
Thus $(P,S')_H=i$ and hence 
$\langle P\rangle_{H[P\cup S']}=S'$ which contradicts 
\eqref{P} and Lemma~\ref{G[X]}. 
It follows that $P\subseteq M$. 

Lastly, it follows from \eqref{I} that $I(G)\subseteq M$. 
Consequently, $M=V(H)$. 
\epr

\begin{cor}\label{Cbound}
For every  graph $G$ such that $|m(G)|\geq 2$, if $\log_2(m(G))$ is not an integer, then 
$p(G)=\lceil\log_2(m(G))\rceil$. 
\end{cor}

\bpr
It suffices to apply Theorem~\ref{mainone} after recalling that 
$\lceil\log_2(m(G))\rceil=\lceil\log_2(m(G)+1)\rceil$ if and only if $\log_2(m(G))$ is not an integer. 
\epr

Before showing Theorem~\ref{maintwo}, we observe 

\begin{lem}\label{isolated}
Given a graph $G$, if ${\rm Iso}(G)\neq\emptyset$ or 
${\rm Iso}(\overline{G})\neq\emptyset$, then 
$$p(G)\geq\lceil\log_2(\max(|{\rm Iso}(G)|,|{\rm Iso}(\overline{G})|)+1)\rceil.$$ 
\end{lem}

\bpr
Let $G$ be a graph such that $\max(|{\rm Iso}(G)|,|{\rm Iso}(\overline{G})|)>0$. 
By interchanging $G$ and $\overline{G}$, assume that ${\rm Iso}(G)\neq\emptyset$. 
Given $p<\lceil\log_2(|{\rm Iso}(G)|+1)\rceil$, 
consider any $p$-extension $H$ of $G$.  
We have $2^{|V(H)\setminus V(G)|}\leq |{\rm Iso}(G)|$ and we verify that $H$ is not prime. 

For each $u\in {\rm Iso}(G)$, we have $N_H(u)\subseteq V(H)\setminus V(G)$. 
Thus $(N_H)_{\restriction {\rm Iso}(G)}$ is a function from ${\rm Iso}(G)$ in 
$2^{V(H)\setminus V(G)}$. 
As previously observed, if $(N_H)_{\restriction {\rm Iso}(G)}$ is not injective, then $\{u,v\}$ is a non-trivial module of $H$ for $u\neq v\in{\rm Iso}(G)$ such that 
$N_H(u)=N_H(v)$. 
So assume that $(N_H)_{\restriction {\rm Iso}(G)}$ is injective. 
As $2^{|V(H)\setminus V(G)|}\leq |{\rm Iso}(G)|$, we obtain that 
$(N_H)_{\restriction {\rm Iso}(G)}$ is bijective. 
Thus there is $u\in{\rm Iso}(G)$ such that $N_H(u)=\emptyset$, that is, $u\in{\rm Iso}(H)$. 
Therefore $H$ is not prime. 
It follows that $p(G)\geq\lceil\log_2(\max(|{\rm Iso}(G)|,|{\rm Iso}(\overline{G})|)+1)\rceil$. 
\epr

We prove Theorem~\ref{maintwo} when $m(G)=2$. 

\begin{prop}\label{m(G)=2}
For every graph $G$ such that $m(G)=2$, 
\begin{equation*}
\text{$p(G)=2$ if and only if $|{\rm Iso}(G)|=2$ or  $|{\rm Iso}(\overline{G})|=2$.}
\end{equation*}
\end{prop}

\bpr 
By Theorem~\ref{mainone}, $p(G)=1$ or 2. 
To begin, assume that $|{\rm Iso}(G)|=2$ or  $|{\rm Iso}(\overline{G})|=2$. 
By Lemma~\ref{isolated}, $p(G)\geq 2$ and hence $p(G)=2$. 
Conversely, assume that $p(G)=2$. 
Let $\alpha\not\in V(G)$. 
As $m(G)=2$, $|C|=2$ for each $C\in\mathcal{C}(G)\cup\mathcal{S}(G)$. 
Let $C_0\in\mathcal{C}(G)\cup\mathcal{S}(G)$. 
We consider any graph $H$ defined on $V(G)\cup\{\alpha\}$ and satisfying the following. 
\begin{enumerate}
\item For each $C\in\mathcal{C}(G)\cup\mathcal{S}(G)$, $\alpha\not\sim_H C$.
\item Let $P\in\mathcal{P}(G)$. We have $G[P]$ is prime. 
Using Lemma~\ref{1extension}, we consider for $H[P\cup\{\alpha\}]$ a prime 1-extension of $G[P]$ to $P\cup\{\alpha\}$.
\item Let $v\in I(G)$. There is $i\in\{0,1\}$ such that $(v,C_0)_G=i$. We require that 
$(v,\alpha)_H\neq i$. 
\item $H[V(G)]=G$.
\end{enumerate}
Since $p(G)=2$, $H$ admits a non-trivial module $M$. 
First, we verify that $\alpha\in M$. 
Otherwise $M$ is a module of $G$. 
By Remark~\ref{base}, there is $C\in\mathcal{C}(G)\cup\mathcal{S}(G)\cup\mathcal{P}(G)$ such that $|C\cap M|\geq 2$. 
Suppose that $C\in\mathcal{C}(G)\cup\mathcal{S}(G)$. 
As $|C|=2$, $C\subseteq M$ which contradicts $\alpha\not\sim_H C$. 
Suppose that $C\in\mathcal{P}(G)$. 
By Lemma~\ref{prime}, $C\cap M=C$. 
Thus $C$ would be a module of $H[C\cup\{\alpha\}]$ which contradicts the fact that 
$H[C\cup\{\alpha\}]$ is prime. 
It follows that $\alpha\in M$. 

Second, we show that $P\subseteq M$ for each $P\in\mathcal{P}(G)$. 
Since $H[P\cup\{\alpha\}]$ is prime and since $M\cap (P\cup\{\alpha\})$ is a module of $H[P\cup\{\alpha\}]$ with 
$\alpha\in M\cap (P\cup\{\alpha\})$, we obtain either $(M\setminus\{\alpha\})\cap P=\emptyset$ or 
$P\subseteq M\setminus\{\alpha\}$. 
In the first instance, there is $i\in\{0,1\}$ such that $(M\setminus\{\alpha\},P)_G=i$. 
Therefore $(\alpha,P)_H=i$ which contradicts the fact that 
$H[P\cup\{\alpha\}]$ is prime. 
It follows that $P\subseteq M$. 

Third, we proof that $C\cap M\neq\emptyset$ for each $C\in\mathcal{C}(G)\cup\mathcal{S}(G)$. 
Otherwise consider $C\in\mathcal{C}(G)\cup\mathcal{S}(G)$ such that $C\cap M=\emptyset$. 
There is $i\in\{0,1\}$ such that $(M,C)_G=i$. 
Thus $(\alpha,C)_H=i$ which contradicts $\alpha\not\sim_H C$. 

In particular we have $C_0\cap M\neq\emptyset$. 
Let $v\in I(G)$. 
Since $(v,C_0\cap M)_G\neq (v,\alpha)_H$, we obtain $v\in M$. 

To conclude, consider $v\in V(H)\setminus M$. 
By what precedes, there is $C\in\mathcal{C}(G)\cup\mathcal{S}(G)$ such that $v\in C$. 
By interchanging $G$ and $\overline{G}$, assume that $C\in\mathcal{S}(G)$. 
Since $v\sim_HM$ and $(v,C\cap M)_G=0$, we obtain $(v,M)_H=0$. 
In particular $(v,I(G))_G=0$ and $(v,P)_G=0$ for every $P\in\mathcal{P}(G)$. 
Let $D\in (\mathcal{C}(G)\cup\mathcal{S}(G))\setminus\{C\}$. 
As $D\cap M\neq\emptyset$, we have $(v,D\cap M)_G=0$ and hence $(v,D)_G=0$ because $D$ is a module of $G$. 
It follows that $v\in{\rm Iso}(G)$. 
Therefore $(C,V(G)\setminus C)_G=0$ because $C$ is a module of $G$. 
Since $C$ is a stable set in $G$, we obtain $C\subseteq {\rm Iso}(G)$. 
Clearly ${\rm Iso}(G)$ is a module of $G$ and a stable set in $G$. 
Thus $|{\rm Iso}(G)|\leq m(G)=2$. 
Consequently $C={\rm Iso}(G)$. 
\epr

We use the following notation in the proof of Theorem~\ref{maintwo}. 
Given a graph $G$ such that $m(G)\geq 3$, set
\begin{align*}
\mathcal{C}_{{\rm max}}(G)=\{C\in\mathcal{C}(G):|C|=m(G)\},\\
\mathcal{S}_{{\rm max}}(G)=\{S\in\mathcal{S}(G):|S|=m(G)\}.
\end{align*}

\bprthtwo
Consider a graph $G$ such that $m(G)=2^k$ where $k\geq 1$. 
By Theorem~\ref{mainone}, $p(G)=k$ or $k+1$. 
To begin, assume that $|{\rm Iso}(G)|=2^k$ or $|{\rm Iso}(\overline{G})|=2^k$. 
By Lemma~\ref{isolated}, $p(G)\geq k+1$ and hence $p(G)=k+1$. 

Conversely, assume that $p(G)=k+1$. 
If $k=1$, then it suffices to apply Proposition~\ref{m(G)=2}. 
So assume that $k\geq 2$. 
With each $C\in\mathcal{C}_{{\rm max}}(G)\cup\mathcal{S}_{{\rm max}}(G)$ associate 
$w_C\in C$. 
Set $W=\{w_C:C\in\mathcal{C}_{{\rm max}}(G)\cup\mathcal{S}_{{\rm max}}(G)\}$ and 
$G'=G[V(G)\setminus W]$. 

We prove that $m(G')=2^k-1$. 
Given $C\in\mathcal{C}_{{\rm max}}(G)\cup\mathcal{S}_{{\rm max}}(G)$, $C\setminus\{w_C\}$ is a module of $G'$ and $C\setminus\{w_C\}$ is a clique or a stable set in $G'$. 
Thus $2^k-1=|C\setminus\{w_C\}|\leq m(G')$. 
Consider $C'\in\mathcal{C}_{{\rm max}}(G')\cup\mathcal{S}_{{\rm max}}(G')$. 
We show that $C'$ is a module of $G$. 
We have to verify that for each $C\in\mathcal{C}_{{\rm max}}(G)\cup\mathcal{S}_{{\rm max}}(G)$, $w_C\sim_G C'$. 
First, asume that there is $c\in(C\setminus\{w_C\})\setminus C'$. 
We have $c\sim_G C'$. 
Furthermore $\{c,w_C\}$ is a module of $G$. 
Thus $w_C\sim_G C'$. 
Second, assume that $C\setminus\{w_C\}\subseteq C'$. 
Clearly $w_C\sim_G C'$ when $C'\subseteq C\setminus\{w_C\}$. 
Otherwise assume that $C'\setminus (C\setminus\{w_C\})\neq\emptyset$. 
By interchanging $G'$ and $\overline{G'}$, assume that $C'$ is a clique in $G'$. 
As $C\setminus\{w_C\}\subseteq C'$ and $|C\setminus\{w_C\}|\geq 2$, we obtain that $C$ is a clique in $G$. 
Since $(C\setminus\{w_C\},C'\setminus C)_G=1$ and since $C$ is a module of $G$, we have 
$(w_C,C'\setminus C)_G=1$. 
Furthermore $(w_C,C\setminus\{w_C\})_G=1$ because $C$ is a clique in $G$. 
Therefore $(w_C,C')_G=1$. 
Consequently $C'$ is a module of $G$. 
As $C'$ is a clique or a stable set in $G$, there is $C\in\mathcal{C}(G)\cup\mathcal{S}(G)$ such that $C\supseteq C'$. 
If $C\not\in\mathcal{C}_{{\rm max}}(G)\cup\mathcal{S}_{{\rm max}}(G)$, then 
$|C'|\leq|C|<m(G)$. 
If $C\in\mathcal{C}_{{\rm max}}(G)\cup\mathcal{S}_{{\rm max}}(G)$, then 
$C'\subseteq C\setminus\{w_C\}$ and hence $|C'|<|C|=m(G)$. 
In both cases, we have $|C'|=m(G')<m(G)$. 
It follows that $m(G')=2^k-1$. 

By Corollary~\ref{Cbound}, $p(G')=k$ and hence there exists a prime $k$-extension $H'$ of $G'$. 
We extend $H'$ to $V(H')\cup W$ as follows. 
Let $C\in\mathcal{C}_{{\rm max}}(G)\cup\mathcal{S}_{{\rm max}}(G)$. 
Consider the function $f_C:C\setminus\{w_C\}\longrightarrow 2^{V(H')\setminus V(G')}$ defined by $c\mapsto N_{H'}(c)\setminus V(G')$ for $c\in C\setminus\{w_C\}$. 
Since $H'$ is prime, $f_C$ is injective. 
As $|C\setminus\{w_C\}|=2^k-1$ and $|2^{V(H')\setminus V(G')}|=2^k$, there is a unique 
$N_C\subseteq V(H')\setminus V(G')$ such that $f_C(c)\neq N_C$ for every 
$c\in C\setminus\{w_C\}$. 
Let $H$ be the extension of $H$ to $V(H')\cup W$ such that 
$N_H(w_C)\cap (V(H')\setminus V(G'))=N_C$ for each 
$C\in\mathcal{C}_{{\rm max}}(G)\cup\mathcal{S}_{{\rm max}}(G)$. 

As $p(G)=k+1$, $H$ is not prime. 
So consider a non-trivial module $M$ of $H$. 
Set $X=V(H')$. We have $H[X]$ is prime. 
Given $u\in X$, we verify that $X_H(u)=\emptyset$. 
Otherwise there is $C\in\mathcal{C}_{{\rm max}}(G)\cup\mathcal{S}_{{\rm max}}(G)$ such that $w_C\in X_H(u)$. 
If $u\in V(G')$, then $\{u,w_C\}$ is a module of $G$. 
Therefore there is $D\in\mathcal{C}(G)\cup\mathcal{S}(G)$ such that $\{u,w_C\}\subseteq D$. 
Necessarily $D=C$ and we would obtain $f_C(u)=N_C$. 
So suppose that $u\in V(H')\setminus V(G')$. 
There is $i\in\{0,1\}$ such that $(w_C,C\setminus\{w_C\})_G=i$. 
Thus $(u,C\setminus\{w_C\})_{H'}=i$. 
Since $f_C$ is injective, the function $g_C:C\setminus\{w_C\}\longrightarrow 2^{((V(H')\setminus V(G'))\setminus\{u\})}$, defined by 
$g_C(c)=f_C(c)\setminus\{u\}$ for $c\in C\setminus\{w_C\}$, is injective as well. 
We would obtain $2^k-1\leq 2^{k-1}$ which does not hold when $k\geq 2$. 
It follows that $X_H(u)=\emptyset$ for each $u\in X$. 
By Lemma~\ref{G[X]}, either $M\cap X=\emptyset$ or $X\subseteq M$. 
In the first instance, $M\subseteq W$ and $M$ is a module of $G$ which contradicts Remark~\ref{base}. 
Consequently $X\subseteq M$. 
As $M$ is a non-trivial module of $H$, there exists $C\in\mathcal{C}_{{\rm max}}(G)\cup\mathcal{S}_{{\rm max}}(G)$ such that $w_C\not\in M$. 
By interchanging $G$ and $\overline{G}$, assume that $C$ is a stable set in $G$. 
We have $(w_C,C\setminus\{w_C\})_G=0$ so that $(w_C,M)_H=0$ and $(w_C,V(G'))_G=0$. 
Given $D\in (\mathcal{C}_{{\rm max}}(G)\cup\mathcal{S}_{{\rm max}}(G))\setminus\{C\}$, we obtain $(w_C,D\setminus\{w_D\})_G=0$. 
Since $D$ is a module of $G$, $(w_C,w_D)_G=0$. 
It follows that $w_C\in{\rm Iso}(G)$. 
As at the end of the proof of Proposition~\ref{m(G)=2}, we conclude by $C={\rm Iso}(G)$. 
\epr

Lastly, we examine the graphs $G$ such that $m(G)=1$. For these, 
$\mathcal{C}(G)=\mathcal{S}(G)=\emptyset$. Thus either $|V(G)|\leq 1$ or $|V(G)|\geq 4$ and $G$ is not prime. 

\begin{prop}\label{m(G)=1}
For every non-prime graph $G$ such that $|V(G)|\geq 4$ and $m(G)=1$, we have $p(G)=1$. 
\end{prop}

\bpr
Since $m(G)=1$, we have $\mathcal{C}(G)=\mathcal{S}(G)=\emptyset$. 
By Corollary~\ref{cpartition}, $\mathfrak{M}(G)=\mathcal{P}(G)\cup\mathcal{I}(G)$. 
By considering $V(G)\in\mathcal{M}(G)$,  it follows from Remark~\ref{base} that there is $P_0\in\mathcal{P}(G)$. 

Let $\alpha\not\in V(G)$. 
We consider any graph $H$ defined on $V(G)\cup\{\alpha\}$ and satisfying the following. 
\begin{enumerate}
\item Let $P\in\mathcal{P}(G)$. 
We have $G[P]$ is prime. 
Using Lemma~\ref{1extension}, we consider for $H[P\cup \{\alpha\}]$ a prime graph such that $H[P]=G[P]$. 

Set $X=P_0\cup \{\alpha\}$. We have $H[X]$ is prime. 
\item Let $v\in I(G)$. Since $P_0$ is a module of $G$ such that $v\not\in P_0$, there is $i\in\{0,1\}$ such that $(v,P_0)_G=i$. We consider for $H[\{v,\alpha\}]$ the graph such that $(v,\alpha)_{H[\{v,\alpha\}]}\neq i$.
\end{enumerate}
In the construction above, we have $H[N]=G[N]$ for each $N\in\mathfrak{M}(G)$. 
Thus we can also assume that $H[V(G)]=G$. 
As in the proof of Theorem~\ref{mainone}, we verify the following. 
For every module $M$ of $H$ such that $|M|\geq 2$, we have 
$\alpha\in M$. 

Now, we prove that $H$ is prime. 
Consider a module $M$ of $H$ such that $|M|\geq 2$. 
By what precedes, $\alpha\in M$ and it follows from Lemma~\ref{G[X]} that either $M\cap X=\{\alpha\}$ or $X\subseteq M$. 
In the first instance, 
there is $i\in\{0,1\}$ such that $(P_0,M\setminus\{\alpha\})_G=1$. 
Thus $(P_0,\alpha)_H=1$ which contradicts the fact that $H[P_0\cup\{\alpha\}]$ is prime. 
It follows that $X\subseteq M$. 
We conclude as in the proof of Theorem~\ref{mainone}. 
Since $(v,P_0)_G\neq (v,\alpha)_H$ for every $v\in I(G)$, we have $I(G)\subseteq M$. 
Lastly, consider $P\in\mathcal{P}(G)\setminus\{P_0\}$. 
As $H[P\cup\{\alpha\}]$ is prime, it follows from Proposition~\ref{properties} that either 
$(M\setminus\{\alpha\})\cap P=\emptyset$ or $P\subseteq M\setminus\{\alpha\}$. 
In the first instance, there is $i\in\{0,1\}$ such that $(P,P_0)_G=i$. 
Therefore $(P,\alpha)_H=1$ which contradicts the fact that $H[P\cup\{\alpha\}]$ is prime. 
Consequently $P\subseteq M$. 
\epr

\end{document}